\documentclass[pdftex,a4paper,11pt]{article}
\usepackage[T1]{fontenc}       
\usepackage{amssymb}
\usepackage{amsmath}
\usepackage{amsthm}
\usepackage{a4wide}
\usepackage{dsfont}            
\usepackage{MnSymbol}          
\usepackage{mathrsfs}          
\usepackage{stmaryrd}          

\newcommand{\R}{\ensuremath{\mathds R}}						
\newcommand{\N}{\ensuremath{\mathds N}}						
\newcommand{\wP}{\ensuremath{\mathds P}}				  

\newcommand{\E}{\ensuremath{\mathds E}}

\newcommand{\one}{\ensuremath{\mathds 1}}         

\newtheorem{thm}{Theorem}

\newtheorem*{errorRep'}{Theorem \ref{errorRep}$'$}
\newtheorem*{result'}{Theorem \ref{result}$'$}
\theoremstyle{definition}
\newtheorem{ex}[thm]{Example}
\newtheorem{rem}[thm]{Remark}

\allowdisplaybreaks

\begin{document}
\date{}
\title{Weak order for the discretization of the stochastic heat equation driven by impulsive noise}
\author{Felix Lindner\thanks{Corresponding author. Tel:  +49-351-463 32425;
fax: +49-351-463 37251.}
  \and René L. Schilling}
\date{\small\textsl{
Institut für Mathematische Stochastik,
Technische Universität Dresden,\\
D-01062 Dresden, Germany}}
\maketitle
\let\thefootnote\relax\footnotetext{\textsl{E-mail addresses:} felix.lindner@tu-dresden.de (F. Lindner), rene.schilling@tu-dresden.de (R. L. Schilling).}
\begin{abstract}
Considering a linear parabolic stochastic partial differential equation driven by impulsive space time noise,
\[dX_t+AX_t\,dt= Q^{1/2}\,dZ_t,\quad X_0=x_0\in H,\quad t\in [0,T],\]
we approximate the distribution of $X_T$.
$(Z_t)_{t\in[0,T]}$ is an impulsive cylindrical process and $Q$ describes the spatial covariance structure of the noise; $\text{Tr}(A^{-\alpha})<\infty$ for some $\alpha>0$ and $A^\beta Q$ is bounded for some $\beta\in(\alpha-1,\alpha]$.

A discretization $(X_h^n)_{n\in\{0,1,\ldots,N\}}$ is defined via the finite element method in space (parameter $h>0$) and a $\theta$-method in time (parameter $\Delta t=T/N$). For $\varphi\in C^2_b(H;\mathds R)$ we show an integral representation for the error $|\mathds E\varphi(X^N_h)-\mathds E\varphi(X_T)|$ and prove that
\[|\mathds E\varphi(X^N_h)-\mathds E\varphi(X_T)|=O(h^{2\gamma}+(\Delta t)^{\gamma})\]
where $\gamma<1-\alpha+\beta$.

\end{abstract}

\bigskip
\noindent
{\small\textsl{2010 MSC: } Primary: 60H15, 65M60; Secondary: 60H35, 60G51, 60G52, 65C30.
\bigskip \\
\noindent
\textsl{Key words and phrases:} Weak order, stochastic heat equation, impulsive cylindrical process, infinite dimensional Lévy process, finite element, Euler scheme.}

\section{Introduction}
In this paper, we study the weak order of convergence of numerical approximations of the solutions of a certain class of linear parabolic stochastic partial differential equations (SPDEs, for short) driven by impulsive space time noise. Unlike the strong order of convergence which measures the pathwise approximation of the true solution by a numerical one (\frenchspacing{cf.}, \frenchspacing{e.g.} \cite{GyMi05}, \cite{GyMi07}, \cite{GyMi09}, \cite{Wal05}, \cite{Yan05}), the weak order is concerned with the approximation of the law of the true solution at a fixed time. There are not many works in literature about the weak approximation of the solutions of SPDEs (see \cite{BouDeb06}, \cite{Deb}, \cite{DebPrin}, \cite{GKL09}, \cite{Haus03}) and, to our knowledge, only SPDEs driven by Gaussian noise have been considered in this context so far.  This work extends the paper $\cite{DebPrin}$ by A. Debussche and J. Printems, where the following Hilbert space valued stochastic differential equation is considered:
\begin{equation}\label{sde0}
dX_t+AX_t\,dt=Q^{1/2}\,dW_t,\quad X_0=x_0\in H,\quad t\in[0,T].
\end{equation}
Here $A:D(A)\subset H\to H$ is a unbounded strictly positive definite self-adjoint operator whose domain $D(A)$ is compactly embedded in $H$; $Q:H\to H$ is a bounded nonnegative definite symmetric operator and $(W_t)_{t\in[0,T]}$ is a cylindrical Wiener process on $H$, $T\in (0,\infty)$. A standard reference for this setting is \cite{DaPraZab}.\\
If we set $H:=L^2(\mathcal O)=L^2(\mathcal O,\mathcal B(\mathcal O),d\xi)$, $\mathcal O\subset \R^d$ open and bounded, and $(A,D(A)):=(-\Delta,\,H^2(\mathcal O)\cap H^1_0(\mathcal O))$, then ($\ref{sde0}$) is an abstract formulation of the stochastic heat equation with Dirichlet boundary conditions
\begin{equation}\label{sde0alt}
\left.
\begin{alignedat}{2}
\frac{\partial X(t,\xi)}{\partial t}-\Delta X(t,\xi)&= \dot\eta(t,\xi), && \qquad (t,\xi)\in[0,T]\times\mathcal O,\\
X(\cdot,\cdot)&= 0 && \qquad \text{ on }\;[0,T]\times\partial\mathcal O,\\
X(0,\cdot)&= x_0 &&\qquad \text{ on }\;\mathcal O.
\end{alignedat}
\quad\right\rbrace
\end{equation}
Here $\dot\eta=\frac{\partial \eta}{\partial t}$ is a generalized random function which can be described as the time derivative (in a distributional sense, cf. \cite{Roz}, \cite{Wal86}) of a real-valued generalized Gaussian process, formally written as
\begin{equation}\label{noise0alt}
\eta(t,\xi)=\int_{\mathcal O}q_0(\xi,\zeta)W(t,\zeta)\,d\zeta,
\end{equation}
where $(W(t,\,\cdot\,))_{t\in[0,T]}$ is a cylindrical Wiener process on $L^2(\mathcal O)$ and $q_0$ is a (generalized) function on $\mathcal O\times \mathcal O$. The operator $Q$ is given by $Qx(\xi)=\int_{\mathcal O}q(\xi,\zeta)x(\zeta)\,d\zeta$ with $q(\xi,\zeta)=\int_{\mathcal O}q_0(\xi,\tau)q_0(\tau,\zeta)\,d\tau$ describing the spatial correlation of the noise $\eta$, cf. \cite{PesZab}, Ch. 4.9.2 and Example 14.26.

Throughout this article, let $H:=L^2(\mathcal O)$. We consider the equation
\begin{equation}\label{sde}
dX_t + AX_t\,dt =  Q^{1/2}\, dZ_t,\quad X_0=x_0\in H,\quad t\in [0,T],
\end{equation}
where $A$ and $Q$ are as above and $(Z_t)_{t\in[0,T]}=(Z(t,\,\cdot\,))_{t\in[0,T]}$ is an impulsive cylindrical process on $H$, see Section 2 for the definition. This is an abstract version of problem ($\ref{sde0alt}$) if one replaces $W(t,\zeta)$ by $Z(t,\zeta)$ in the formal definition ($\ref{noise0alt}$) of the noise $\eta$ and if furthermore $q_0$ is symmetric and positive semidefinite in the sense that $\int_{\mathcal O}\int_{\mathcal O}q_0(\xi,\zeta)\psi(\xi)\psi(\zeta)\,d\xi\,d\zeta \geq 0$ for all test functions $\psi\in C_0^\infty(\mathcal O)$.

In $\cite{DebPrin}$, a discretization $(X_h^n)_{n\in\{1,\ldots,N\}}$ of the solution $(X_t)_{t\in[0,T]}$ of equation ($\ref{sde0}$) is obtained by the finite element method in space (parameter $h>0$) and a $\theta$-method in time (parameter $\Delta t=T/N$). Under the assumption that $A^{-\alpha}$ is a finite trace operator for some $\alpha>0$ and that $A^\beta Q$ is bounded for some $\beta\in (\alpha-1,\alpha]$, it is shown there that for functions $\varphi\in C_b^2(H,\R)$,
\begin{equation}\label{order}
|\E\varphi(X^N_h)-\E\varphi(X_T)|\leq C\cdot (h^{2\gamma}+(\Delta t)^\gamma)
\end{equation}
for any $\gamma<1-\alpha+\beta\leq 1$.
\bigskip

In this paper, we consider the analoguous discretization of the solution of ($\ref{sde}$) and make the same assumptions on the operators $A$ and $Q$ as in $\cite{DebPrin}$. We give a representation formula for the error (Theorem~\ref{errorRep}) and, under some integrability condition on the jump size intensity $\nu$ of the cylindrical impulsive process $(Z_t)_{t\in [0,T]}$, we show that ($\ref{order}$) holds also for the solution $(X_t)_{t\in[0,T]}$ in ($\ref{sde}$) and the corresponding discretization (Theorem~\ref{result}).

SPDEs driven by impulsive noise (or Poisson noise) have been considered, {\frenchspacing e.g. }in \cite{AppWu}, \cite{Haus05}, \cite{Knoche}, \cite{MPR08}, \cite{Mue98}, \cite{Myt02}. The monograph \cite{PesZab} gives a good overview about SPDEs driven by Lévy noise. In \cite{Haus06} and \cite{Haus08}, numerical approximations in time and space of SPDEs driven by Poisson random measures are investigated and the strong error is estimated. Of course, this especially implies an estimate for the weak approximation error. A difference to our result is that we look at impulsive noise which is white in time and coloured in space whereas in \cite{Haus06} and \cite{Haus08} a class of SPDEs driven by Poisson random measures which correspond to impulsive space time white noise is considered. Our motivation for this paper was to show that the techniques applied in \cite{DebPrin} with respect to the cylindrical Wiener process also work for certain jump processes.

The main technical difference between $(\ref{sde})$ and $(\ref{sde0})$ lies in the fact that the impulsive cylindrical process $(Z_t)_{t\in[0,T]}$ is a purely discontinuous Hilbert space valued martingale, while the cylindrical Wiener process $(W_t)_{t\in[0,T]}$ is continuous. As a consequence, the main tools for estimating the weak order of convergence for the numerical scheme --- the Itô formula and (connected with it) the backward Kolmogorov equations for certain processes associated with the solutions of the SPDEs and their discretizations --- are completely different for $(\ref{sde})$ and $(\ref{sde0})$. The main task therefore is to find manageable expression for the approximation error, which allows estimates using techniques similar to those in $\cite{DebPrin}$.

We note that ($\ref{order}$) remains true for the solution $(X_t)_{t\in[0,T]}$ of
\begin{equation*}\label{sdeComb}
dX_t + AX_t\,dt =  Q^{1/2}_0\,dW_t + Q_1^{1/2}\, dZ_t,\quad X_0=x_0\in H,\quad t\in [0,T],
\end{equation*}
and the corresponding discretization. Here, the covariance operators $Q_0$ and $Q_1$ are assumed to have the same properties as $Q$ above.

Let us finally remark that the weak order for numerical approximations of ordinary stochastic differential equations with jumps has been studied in several papers, see \cite{BruPla} for a comprehensive survey. In \cite{ProTal} the Euler scheme for Lévy driven equations of the form $dX_t=f(X_{t-})\,dZ_t$ is considered and it is proved that the weak order is $O(\Delta t)$ if the jump intensity measure of $Z$ has its first several moments finite. In \cite{MicPla} and \cite{LiuLi}, the more general concept of stochastic Taylor expansions is used but the underlying jump intensity measure is supposed to be finite,  \frenchspacing{i.e. only finitely} many jumps appear on a bounded time interval. Comparing these finite dimensional settings with \eqref{sde}, we note the following: If $H$ is finite dimensional the conditions on $A$ and $Q$ will be trivially fulfilled for $\alpha=\beta=1$ and one will get the same order as in \cite{ProTal}. Secondly, the jump intensity measure of the cylindrical process $(Z_t)_{t\in[0,T]}$ with state space $U\supset H=L^2(\mathcal O)$ is a measure on $(U,\mathcal B(U))$ and has to be distinguished from the jump size intensity measure $\nu$ as a measure on $(\R,\mathcal B(\R))$, see Section 2. The conditions we impose on $(Z_t)_{t\in[0,T]}$ in order to obtain our results can be best expressed in terms of moments of the jump size intensity $\nu$. However, the jump intensity measure and the jump size intensity measure of the impulsive cylindrical process are closely connected, see Remark \ref{jumpIntMeas}. In particular, the jump intensity measure of the driving process $(Z_t)_{t\in[0,T]}$ in \eqref{sde} is allowed to be infinite. Thirdly, having found a suitable error expansion, the following estimates contain additional difficulties due to the infinite-dimensionality of our problem. This is why in Section 6 the assumptions on the jump intensity measure of the driving process are stronger than those in \cite{ProTal}.

\section{Impulsive cylindrical process}
Let $\mathcal O\subset\R^d$ be open and bounded, $T\in (0,\infty)$, and consider the product space $[0,T]\times\mathcal O\times\R$ equipped with the Borel $\sigma$-algebra $\mathcal B([0,T]\times\mathcal O\times \R)$. The generic element in $[0,\infty)\times\mathcal O\times\R$ is denoted by $(t,\xi,\sigma)$ or $(s,\xi,\sigma)$. Let $\nu$ be a sigma-finite measure on $\R$ and $\pi$ a Poisson random measure on  $[0,\infty)\times\mathcal O\times\R$ with reference measure $dt\,d\xi\,\nu(d\sigma)$. By $\hat\pi$ we denote the compensated Poisson random measure, i.e. \[\hat\pi(dt,d\xi,d\sigma)=\pi(dt,d\xi,d\sigma)-dt\,d\xi\,\nu(d\sigma).\]
Let $\pi$ be defined on a complete probability space $(\Omega,\mathcal A,\wP)$ with a filtration $(\mathcal F_t)_{t\in[0,T]}$ satisfying the usual hypotheses, cf. \cite{Met}. We assume that $\pi([0,t]\times B)$ is $(\mathcal F_t)$-measurable for all $t\in[0,T],\;B\in\mathcal B(\mathcal O\times \R)$ and that $\pi((s,t]\times B)$ is independent of $\mathcal F_s$ for $0\leq s<t\leq T,\;B\in\mathcal B(\mathcal O\times \R)$. Furthermore, we assume that $\pi$ is of the form
\begin{equation}
\pi(A)(\omega)=\sum_{j=1}^\infty\delta_{(T_j(\omega),\Xi_j(\omega),\Sigma_j(\omega))}(A),\qquad A\in\mathcal B([0,T]\times\mathcal O\times\R),\;\omega\in\Omega,
\end{equation}
for a properly chosen sequence $\big((T_j,\Xi_j,\Sigma_j)\big)_{j\in\N}$ of random elements in $[0,T]\times\mathcal O\times\R$, cf. \cite{PesZab}, Chapter 6.

For general impulsive cylindrical processes, $\nu$ is a Lévy measure, i.e. $\nu(\{0\})=0$ and $\int_\R\min\big(\sigma^2, 1\big)\,\nu(d\sigma)<\infty$, cf. \cite{PesZab}, Ch. 7.2. To simplify the exposition we will additionally assume that $\int_{|\sigma|\geq 1}\sigma^2\,\nu(d\sigma)<\infty$ which is equivalent to saying that
\begin{equation}\label{AssNu0}
\int_\R\sigma^2\,\nu(d\sigma)<\infty.
\end{equation}
In the Appendix \ref{A2}, we show how our results can be extended to the general case.
Under condition \eqref{AssNu0} the random variables $Z_t^{(k)}$ and $Z_t$ defined below have finite second moments. (This is, {\frenchspacing e.g. }the case if $\nu$ is a Lévy measure with bounded support.)

To fix notation, let us give a brief survey of $L^2$-integration {\frenchspacing w.r.t. $\hat\pi$} for deterministic integrands. Let $f:[0,T]\times\mathcal O\times\R\to\R$ be a simple function, i.e.
\[f=\sum_{k=1}^na_k\one_{A_k},\]
where $n\in\N$, $a_k\in\R$ and $A_k\in\mathcal B([0,T]\times\mathcal O\times\R)$, such that
\[\int_0^T\int_{\mathcal O}\int_\R\one_{A_k}(t,\xi,\sigma)\,\nu(d\sigma)\,d\xi\,dt<\infty,\]
for $k=1,\ldots,n$.\\
For simple functions $f$ the stochastic integral {\frenchspacing w.r.t. $\hat\pi$} is defined by
\[\int_0^T\int_{\mathcal O\times\R}f(t,\xi,\sigma)\,\hat\pi(dt,d\xi,d\sigma):=
\sum_{k=1}^na_k\hat\pi(A_k).\]
Since
\[\E\left|\int_0^T\int_{\mathcal O\times\R}f(t,\xi,\sigma)\,\hat\pi(dt,d\xi,d\sigma)\right|^2=
\int_0^T\int_{\mathcal O\times\R}|f(t,\xi,\sigma)|^2\,dt\,d\xi\,\nu(d\sigma),\]
the stochastic integral can be uniquely extended to an isometric linear operator mapping $L^2([0,T]\times\mathcal O\times\R,ds\,d\xi\,\nu(d\sigma);\R)$ to $L^2(\Omega,\mathcal A,\wP;\R)$. Later, we will also consider $\hat\pi$-integrals of $H$-valued square integrable functions which are defined similarly. In the course of the proof of our result, we will have to deal with $L^1$-integrals of stochastic integrands against the (not compensated) random measure $\pi$; see \cite{PesZab}, Sections 6.2, 8.7 and \cite{Wal86}, Chapter 2, for a detailed exposition of the various integrals.\\
As usual, for $t\in[0,T]$ and $f\in L^2([0,T]\times\mathcal O\times\R,ds\,d\xi\,\nu(d\sigma);\R)$, we define
\[\int_0^t\int_{\mathcal O\times\R}f(s,\xi,\sigma)\,\hat\pi(ds,d\xi,d\sigma):=
\int_0^T\int_{\mathcal O\times\R}\one_{[0,t]}(s)f(s,\xi,\sigma)\,\hat\pi(ds,d\xi,d\sigma).\]

Now we can to define the impulsive cylindrical process. Let $(e_k)_{k\in\N}$ be an orthonormal basis of $H$ and for $k\in\N,\;t\in[0,T]$, let
\begin{equation*}
Z^{(k)}_t:=Z(t,e_k):=\int_0^t\int_{\mathcal O\times\R}e_k(\xi)\cdot\sigma\,\hat\pi(ds,d\xi,d\sigma).
\end{equation*}
Note that the processes $(Z^{(k)}_t)_{t\in[0,T]},\;k\in\N,$ are real-valued, square integrable Lévy processes which are also martingales.\\
Now let $U$ be a further Hilbert space such that $H$ is densely embedded in $U$ and such that the embedding is Hilbert-Schmidt, e.g. $U=H^{-\frac d2-\epsilon}(\mathcal O)$ for some $\epsilon>0$. In \cite{PesZab}, Chapter 7, it is shown that
\begin{equation*}
Z_t:=L^2(\Omega,\mathcal A,\wP;U)\text{-}\lim_{n\uparrow\infty}\sum_{k=1}^n Z^{(k)}_te_k,\qquad t\in\,[0,T],
\end{equation*}
defines an $U$-valued $L^2(\wP)$-Lévy-martingale with reproducing kernel Hilbert space $\big(\mathcal H,\langle\,\cdot\,,\,\cdot\,\rangle_{\mathcal H}\big)=\big(H,\int_\R\sigma^2\,\nu(d\sigma)\langle\,\cdot\,,\,\cdot\,\rangle_H\big)$, $H=L^2(\mathcal O)$. $(Z_t)_{t\in[0,T]}$ is called \emph{impulsive cylindrical process on $L^2(\mathcal O)$ with jump size intensity $\nu$}.

\begin{rem}\label{jumpIntMeas}
The jump intensity measure $\mu$ of the $U$-valued Lévy-process $(Z_t)_{t\in[0,T]}$ is given by
\[
\mu(B)=\lambda\otimes\nu\left(\left\{(\xi,\sigma)\in\mathcal O\times\R\,:\,\sum_{k=1}^\infty e_k(\xi)\sigma e_k\;\in\,B\right\}\right),\quad B\in\mathcal B(U)\]
where $\lambda$ denotes Lebesgue measure on $\mathcal O$ and the infinite sum is a limit in $L^2(\mathcal O\times\R,\,d\xi\,\nu(d\sigma);\,U)$.
\end{rem}

As for the following examples, compare \cite{Mue98}, \cite{Myt02}.

\begin{ex}
Let $\beta\in (0,2),\;\tau\in(0,\infty)$ and consider the jump size intensity
\[\nu(d\sigma)=\frac{1}{\sigma^{1+\beta}}\one_{[0,\tau]}(\sigma)\,d\sigma.\]
Then, for every $B\in \mathcal B(\mathcal O)$, the one-dimensional process \[\big(Z(t,\one_B)\big)_{t\in[0,T]}=\left(\int_0^t\int_{\mathcal O\times\R}\one_B(\xi)\cdot\sigma\,\hat\pi(ds,d\xi,d\sigma)\right)_{t\in[0,T]}\] can be characterized in terms of one-sided $\beta$-stable Lévy processes as follows:\\
If $\beta\in(1,2)$, then $\big(Z(t,\one_B)\big)_{t\in[0,T]}$ is equivalent to the process obtained by removing all jumps greater than $\tau$ from a one-sided $\beta$-stable Lévy process $(L_t)_{t\in[0,T]}$ with Laplace transform
\[\E e^{-rL_t}=\exp\left\{-t|B|c_\beta\cdot r^{\beta}\right\}=\exp\left\{-t|B|\int_0^\infty\big(1-e^{-r\sigma}-r\sigma\big)
\frac{d\sigma}{\sigma^{1+\beta}}\right\},\qquad\;r>0,\]
and then adding the shift $t\mapsto t|B|\int_\tau^\infty \sigma\frac{d\sigma}{\sigma^{1+\beta}}=t|B|\frac{\tau^{1-\beta}}{\beta-1}$. Here $|B|$ denotes the Lebesgue measure of $B\in\mathcal B(\mathcal O)$.\\
If $\beta\in(0,1)$, $\big(Z(t,\one_B)\big)_{t\in[0,T]}$
is equivalent to the process obtained by removing all jumps greater than $\tau$ from a positive one-sided $\beta$-stable Lévy process $(L_t)_{t\in[0,T]}$ with Laplace transform
\[\E e^{-rL_t}=\exp\left\{-t|B|c_\beta\cdot r^{\beta}\right\}=\exp\left\{-t|B|\int_0^\infty\big(1-e^{-r\sigma}\big)
\frac{d\sigma}{\sigma^{1+\beta}}\right\},\qquad\;r>0,\]
and then subtracting the shift $t\mapsto \E\int_0^t\int_{\mathcal O}\int_\R\one_B(\xi)\cdot\sigma\,\pi(ds,d\xi,d\sigma)= t|B|\frac{\tau^{1-\beta}}{1-\beta}$.
\end{ex}

\begin{ex}
Let $\alpha\in(0,2),\;\tau\in(0,\infty)$ and set
\[\nu(d\sigma)=\frac{1}{|\sigma|^{1+\alpha}}\one_{[-\tau,\tau]}(\sigma)\,d\sigma.\]
Then, for $B\in\mathcal B(\mathcal O)$, we have
\begin{equation}\label{zeroComp}
\begin{aligned}
Z(t,\one_B)&=\int_0^t\int_{\mathcal O\times\R}\one_B(\xi)\cdot\sigma\,\hat\pi(ds,d\xi,d\sigma)\\
&=L^2(\Omega,\mathcal A,\wP;\,\R)\text{-}\lim_{\epsilon\searrow 0}\int_0^t
\int_{\mathcal O\times\{|\sigma|\geq\epsilon\}}\one_B(\xi)\cdot\sigma\,\pi(ds,d\xi,d\sigma),
\end{aligned}
\end{equation}
and $\big(Z(t,\one_B)\big)_{t\in[0,T]}$ is equivalent to the process obtained by removing all jumps of absolute value greater than $\tau$ from a symmetric $\alpha$-stable Lévy process $(L_t)_{t\in[0,T]}$ with Fourier transform
\[\E e^{irL_t}=\exp\left\{-t|B|C_\alpha\cdot r^\alpha\right\}=\exp\left\{-t|B|\int_\R(1-\cos(r\sigma))\frac{d\sigma}{|\sigma|^{1+\alpha}}\right\},
\qquad r\in\R,\]
i.e.
\[\E e^{irZ(t,\one_B)}=\exp\left\{-t|B|\int_{-\tau}^\tau(1-\cos(r\sigma))\frac{d\sigma}{|\sigma|^{1+\alpha}}\right\},
\qquad r\in\R.\]
The second equality in \eqref{zeroComp} holds since $\nu$ is symmetric and therefore \[\int_0^t
\int_{\mathcal O\times\{|\sigma|\geq\epsilon\}}\one_B(\xi)\cdot\sigma\;ds\,d\xi\,\nu(d\sigma)=0,\qquad \epsilon>0,\]
so that the integrals w.r.t. $\hat\pi$ and $\pi$ coincide for integrands in $L^2([0,T]\times\mathcal O\times\R,\,ds\,d\xi\,\nu(d\sigma);\,\R)\cap L^1([0,T]\times\mathcal O\times\R,\,ds\,d\xi\,\nu(d\sigma);\,\R)$.
\end{ex}
\bigskip

For the estimate of the weak order of approximation (Theorem~\ref{result}), we have to make a further restriction on the jump size intensity $\nu$ by assuming that
\begin{equation}\label{AssNu}
\int_\R\max\big(|\sigma|,\sigma^2\big)\,\nu(d\sigma)<\infty.
\end{equation}
This is, {\frenchspacing e.g. }fulfilled if $\beta\in(0,1)$ and $\alpha\in(0,1)$ in Examples 1 and 2.

\section{Notation, assumptions and preliminary results}
The inner product and norm of $H=L^2(\mathcal O)$ are denoted by $\langle\cdot,\cdot\rangle_H$ and $|\,\cdot\,|_H$, respectively; $Q:H\to H$ is a bounded nonnegative definite symmetric operator and $A:D(A)\subset H\to H$ a (unbounded) strictly positive definite self-adjoint operator whose domain $D(A)$ (endowed with the graph norm $|\,\cdot\,|_H+|A\,\cdot\,|_H$) is compactly embedded in $H$. Therefore, the spectrum consists of a sequence $(\lambda_k)_{k\in\N}\subset (0,\infty)$. We assume that the eigenvalues are ordered increasingly, $\lambda_1\leq\lambda_2\leq\ldots$, including multiplicities. By $(\tilde e_k)_{k\in\N}$ we denote the corresponding orthonormal basis of eigenvectors.
For any $s\geq 0$ we set
\begin{eqnarray*}
D(A^s)&:=&\left\{u=\sum_{k=1}^\infty \langle u,\tilde e_k\rangle_H\tilde e_k\in H \,:\;\sum_{k=1}^\infty\lambda_k^{2s}\langle u,\tilde e_k\rangle_H^2<\infty\right\},\\
A^su&:=&\sum_{k=1}^\infty\lambda^s_k\langle u,\tilde e_k\rangle_H\tilde e_k,\;\quad \,u\in D(A^s),
\end{eqnarray*}
so that $D(A^s)$ endowed with the graph norm $|\,\cdot\,|_H+|A^s\,\cdot\,|_H$ is a Hilbert space. Furthermore we define $D(A^{-s})$ for $s\geq 0$ as the completion of $H$ with respect to the norm $|\,\cdot\,|_H+|A^{-s}\,\cdot\,|_H$, defined on $H$ by $|A^{-s}u|^2_H= \sum_{k=1}^\infty\lambda_k^{-2s}\langle u,\tilde e_k\rangle_H^2$.

Let $(V_h)_{h>0}$ be a family of finite dimensional subspaces of $V:=D(A^{1/2})$ parameterized by a small parameter $h>0$, given by the finite element method. As standard references for the finite element method we mention \cite{Ciar} and \cite{StrangFix}. By $P_h$ we denote the orthogonal projector from $H$ onto $V_h$ with respect to the inner product $\langle\,\cdot\,,\,\cdot\,\rangle_H$ and $\Pi_h$ is the orthogonal projector from $H$ onto $V_h$ with respect to the inner product $\langle A^{1/2}\,\cdot\,,A^{1/2}\,\cdot\,\rangle_H$. For any $h>0$, we define a positive symmetric bounded linear operator $A_h:V_h\to V_h$ by
\begin{equation*}
\langle A_hu_h,v_h\rangle_H=\langle A^{1/2}u_h,A^{1/2}v_h\rangle_H\qquad\forall (u_h,v_h)\in V_h\times V_h.
\end{equation*}
As in \cite{DebPrin}, let $(S(t))_{t\geq 0}$ denote the $C_0$-semigroup of contractive operators on H generated by $-A$ and let $(S_h(t))_{t\geq 0}$ denote the semigroup of operators on $V_h$ generated by $-A_h$.

We assume that the finite dimensional spaces $V_h$ admit the following two estimates, which are classical properties for standard finite element spaces.
\bigskip\\
\begin{samepage}
{\bf Assumption:} \\
For all $q\in[0,2]$ there exist constants $\kappa_1>0,\;\kappa_2>0$ independent from $h$ such that for all $t>0$
\begin{eqnarray}
\|S_h(t)P_h-S(t)\|_{L(H)}&\leq&\kappa_1 h^qt^{-q/2},\label{S1}\\
\|S_h(t)P_h-S(t)\|_{L(H,D(A^{1/2}))}&\leq&\kappa_2 ht^{-1}.\label{S2}
\end{eqnarray}
\end{samepage}
In \cite{Thom}, Theorem 3.5, \cite{Bram77}, Theorem 3.2, and \cite{John}, Theorem 4.1, the above estimates \eqref{S1} and \eqref{S2} are shown for the case $(A,D(A))=(-\Delta,\,H^2(\mathcal O)\cap H^1_0(\mathcal O))$ under the assumption that
\begin{eqnarray}
 |\Pi_hv-v|_H&\leq&\kappa_0h^s|A^{s/2}v|_H,\label{Vh1}\\
 |A^{1/2}(\Pi_hv-v)|_H&\leq&\kappa_0h^{s-1}|A^{s/2}v|_H\label{Vh2}
\end{eqnarray}
for all  $s\in [1,2],\;v\in D(A^{s/2})$ and some constant $\kappa_0>0$.
Finite elements satisfying \eqref{Vh1} and \eqref{Vh2} are well known, like $P_k$ triangular finite elements on a convex polygonal domain or $Q_k$ rectangular finite elements on a rectangular domain, $k\geq1$, see \cite{Ciar}, \cite{StrangFix}.

The main assumptions concerning the operator $A$ and the covariance operator $Q$ of the noise are the following: There exist real numbers
\begin{equation}
\alpha>0,\quad\beta\in(\alpha-1,\alpha]\label{ass1}
\end{equation}
such that
\begin{gather}
\text{Tr}(A^{-\alpha})=\sum_{n=1}^\infty\lambda_n^{-\alpha}<\infty\;,\label{assTr}\\
A^\beta Q\in L(H).\label{assQ}
\end{gather}
Notice that (\ref{assQ}) implies for any $\lambda\in[0,1]$
\begin{equation}\label{interpol}
A^{\lambda\beta}Q^\lambda\in L(H)\;\text{ and } \;\|A^{\lambda\beta}Q^\lambda\|_{L(H)}\leq\|A^\beta Q\|_{L(H)}^\lambda.
\end{equation}
According to \cite{PesZab}, Chapter 9, equation (\ref{sde}) has a unique (predictable) \emph{weak solution} which is given as the \emph{mild solution}
\begin{equation}
X_t=S(t)x_0 + \int_0^tS(t-s)Q^{1/2}\,dZ_s,\qquad t\in[0,T],
\end{equation}
provided that
\begin{equation}\label{intCond}
\|S(t)Q^{1/2}\|_{\text{(HS)}}\in L^2([0,T],dt),
\end{equation}
where $\|\,\cdot\,\|_{\text{(HS)}}$ is the Hilbert-Schmidt norm.
It is shown in \cite{DebPrin} that \eqref{ass1}, (\ref{assTr}) and (\ref{assQ}) are sufficient conditions for (\ref{intCond}). Throughout this paper, we only consider weak respectively mild solutions of SPDEs.

We proceed with some further notation: By $C_b^k(H)=C^k_b(H;\R)$ we denote the space of all $k$-times continuously Fréchet-differentiable real valued functions on $H$ which are bounded together with their derivatives. For $\phi\in C^1_b(H)$ and $x\in H$ the first-order derivative $D\phi(x)$ of $\phi$ in $x$ is identified with its gradient and thus considered as an element in $H$. Similarly, for $\phi\in C^2_b(H)$ and $x\in H$, the second-order derivative $D^2\varphi(x)$ is seen as an element in $L(H)=L(H,H)$, the space of all linear and bounded operators on $H$.\\
Let $L_{\text{(HS)}}(H)=L_{\text{(HS)}}(H,H)$ denote the Hilbert space of all Hilbert-Schmidt operators on $H$. It is a subspace of $L(H)$. Given an orthonormal basis $(e_k)_{k\in\N}$ of $H$, the scalar product in $L_{\text{(HS)}}(H)$ of operators $T\in L_{\text{(HS)}}(H),\;S\in L_{\text{(HS)}}(H)$ is given by
\[\langle T,S\rangle_{\text{(HS)}}:=\sum_{k\in\N}\langle Se_k,Te_k\rangle_H.\]
The corresponding Hilbert-Schmidt norm is denoted by $\|\;\cdot\;\|_{\text{(HS)}}$.

Finally, we use $\mathcal M^2_T(H)$ for the space of all right continuous $L^2(\wP)$-martingales $M=(M_t)_{t\in[0,T]}$ with values in $H$. $(\llbracket M \rrbracket_t)_{t\in[0,T]}$ denotes the tensor quadratic variation (or operator square bracket) of $M\in\mathcal M^2_T(H)$; $(\llbracket M \rrbracket^c_t)_{t\in[0,T]}$ denotes the continuous part. Note that for $\wP$-almost every $\omega\in\Omega$, for every $t\in [0,T]$, $\llbracket M \rrbracket_t(\omega)$ and $\llbracket M \rrbracket^c_t(\omega)$ belong to the space of nuclear operators on $H$, a subspace of $L_{\text{(HS)}}(H)$ which is continuously embedded in $L_{\text{(HS)}}(H)$, cf. \cite{Met}.

\section{Approximation scheme}\label{ApproximationScheme}
In order to approximate the mild solution
\begin{equation*}
X_t = S(t)x_0 + \int_0^tS(t-s)Q^{1/2}\,dZ_s,\qquad t\in[0,T],
\end{equation*}
of equation $(\ref{sde})$, we adapt the numerical scheme described in $\cite{DebPrin}$, with $(W_t)_{t\in[0,T]}$ replaced by the jump process $(Z_t)_{t\in[0,T]}$.\\
Given an integer $N\geq 1$, set $\Delta t=T/N$ and $t_n=n\Delta t,\;n=0,\ldots,N$. For any $h>0$, the approximations $X_h^n$ of $X_{t_n}$ in $V_h$, $n=0,\ldots,N$, are defined as those $V_h$-valued random variables which satisfy for all $v_h \in V_h$
\begin{align}
\left\langle X_h^{n+1}-X_h^n\,,\,v_h\right\rangle_H+\Delta t\left\langle A^{1/2}\big(\theta X_h^{n+1}+(1-\theta)X_h^n\big)\,,\,A^{1/2} v_h\right\rangle_H
&=\left\langle Q^{1/2}Z_{t_{n+1}}-Q^{1/2}Z_{t_n}\,,\,v_h\right\rangle_H\label{scheme}
\end{align}
and
\begin{equation}\langle X_h^0,v_h\rangle_H=\langle x_0,v_h\rangle_H,\label{initial}\end{equation}
with
\begin{equation}\label{theta}
\theta\in\;(1/2,1].
\end{equation}
Here the expression $\left\langle Q^{1/2}Z_{t_{n+1}}-Q^{1/2}Z_{t_n}\,,\,v_h\right\rangle_H$ is the usual shorthand for
\begin{equation}\label{SP}
L^2(\Omega,\mathcal A,\wP;\,\R)\text{-}\lim_{M\to\infty}\sum_{k=1}^M\big(Z(t_{n+1},e_k)-Z(t_n,e_k)\big)\cdot\left\langle Q^{1/2}e_k\,,\,v_h\right\rangle_H,
\end{equation}
$(e_k)_{k\in\N}$ being an arbitrary orthonormal basis of $H$. In particular, one can choose $(e_k)_{k\in\N}$ as an orthonormal basis consisting solely of eigenvectors of $Q$ and of elements of the kernel of $Q$.\\
Note that $(X_h^n)_{n\in\{0,\ldots,N\}}$ is a discretization of $(X_t)_{t\in[0,T]}$, both in time and space.
Solving the equations above leads to
\begin{eqnarray}
X_h^n &=& S_{h,\Delta t}^n P_hx_0 +\sum_{k=0}^{n-1}S_{h,\Delta t}^{n-k-1}T_{h,\Delta t}P_h Q^{1/2}(Z_{t_{k+1}}-Z_{t_k}),\label{discrSoln}
\end{eqnarray}
for $n\in\{0,\ldots,N\}$ and $h>0$, where
\begin{eqnarray*}
S_{h,\Delta t} &:=&(I+\theta\Delta t A_h)^{-1}(I-(1-\theta)\Delta tA_h),\\
T_{h,\Delta t}&:=&(I+\theta\Delta t A_h)^{-1}.
\end{eqnarray*}
\begin{rem}
In this paper, we are interested in the weak order of the scheme \eqref{scheme}-\eqref{initial}. Theorem~\ref{result} below gives an estimate for the difference $|\E\varphi(X_h^N)-\E\varphi(X_T)|$ for suitable real valued functions $\varphi$.
Obviously, for the numerical implementation one has to truncate the infinite sum in \eqref{SP}. That is, one has to replace the right hand side of \eqref{scheme} by
\[\left\langle \big(Q^{(M)}\big)^{1/2}Z_{t_{n+1}}-\big(Q^{(M)}\big)^{1/2}Z_{t_n}\,,\,v_h\right\rangle_H=\sum_{k=1}^M\big(Z(t_{n+1},e_k)-Z(t_n,e_k)\big)
\cdot\left\langle Q^{1/2}e_k\,,\,v_h\right\rangle_H,\]
for some $M\in\N$.
Here, the operator $Q^{(M)}$ is defined by
\[Q^{(M)}u=\sum_{k=1}^M\langle u,e_k \rangle_HQe_k,\quad u=\sum_{k=1}^\infty \langle u,e_k \rangle_He_k\,\in\,H.\]
Let $\big(X_h^{n,(M)}\big)_{n\in\{0,\ldots,N\}}$ be the solution of the truncated scheme, i.e.
\[X_h^{n,(M)} = S_{h,\Delta t}^n P_hx_0 +\sum_{k=0}^{n-1}S_{h,\Delta t}^{n-k-1}T_{h,\Delta t}P_h \big(Q^{(M)}\big)^{1/2}(Z_{t_{k+1}}-Z_{t_k}),\]
and let
\begin{equation*}
X_t^{(M)} = S(t)x_0 + \int_0^tS(t-s)\big(Q^{(M)}\big)^{1/2}\,dZ_s,\qquad t\in[0,T].
\end{equation*}
Then, in the setting of Theorem~\ref{result} below, it is not hard to see that
\[\left|\E\varphi\big(X_h^{N,(M)}\big)-\E\varphi\big(X_T^{(M)}\big)\right|\leq
 C\cdot (h^{2\gamma}+(\Delta t)^{\gamma}),\]
where $C>0$ and $\gamma>0$ are the same numbers that appear in the upper bound for the error
$|\E\varphi(X_{h}^N)-\E\varphi(X_T\big)|$ in Theorem~\ref{result}. In particular, the constant $C$ does not depend on $M\in \N$.
Therefore, the error $\big|\E\varphi\big(X_h^{N,(M)}\big)-\E\varphi\big(X_T\big)\big|$, which is relevant for numerical simulations, can be estimated by
\[\left|\E\varphi\big(X_h^{N,(M)}\big)-\E\varphi\big(X_T\big)\right|\leq
 C\cdot \big(h^{2\gamma}+(\Delta t)^{\gamma}\big)+\left|\E\varphi\big(X_T^{(M)}\big)
-\E\varphi\big(X_T\big)\right|.
\]
\end{rem}\bigskip

It will be useful to introduce some further notation. For each $h>0$, we will also consider the following spatial discretization of $(X_t)_{t\in[0,T]}$
\begin{eqnarray}\label{Xht}
X_{h,t} &=& S_h(t)P_hx_0 + \int_0^tS_h(t-s)P_hQ^{1/2}\,dZ_s,\qquad t\in [0,T],
\end{eqnarray}
and  the two auxiliary processes
\begin{eqnarray}
Y_{h,t} &=& S_h(T)P_hx_0+\int_0^tS_h(T-s)P_hQ^{1/2}\,dZ_s\notag\\
&=& S_h(T)P_hx_0+\int_0^t\Phi(s)\,dZ_s,\qquad\qquad\qquad\qquad\; t\in [0,T],\label{Yht}\\
\bar{Y}_{h,t}&=& S^N_{h,\Delta t}P_hx_0 + \int_0^t\sum_{k=0}^{N-1}S_{h,\Delta t}^{N-k-1}T_{h,\Delta t}\one_{(t_k,t_{k+1}]}(s)P_hQ^{1/2}dZ_s\notag\\
&=& S^N_{h,\Delta t}P_hx_0+\int_0^t\Gamma(s)\,dZ_s,\qquad\qquad\qquad\qquad\quad t\in [0,T]\label{Ybarht}.
\end{eqnarray}
To ease notation, we have used
\begin{eqnarray*}
\Phi(s)&:=&S_h(T-s)P_hQ^{1/2},\\
\Gamma(s)&:=&\sum_{k=0}^{N-1}S_{h,\Delta t}^{N-k-1}T_{h,\Delta t}\one_{(t_k,t_{k+1}]}(s)P_hQ^{1/2}.
\end{eqnarray*}
Moreover, it will be convenient to set
\[\tilde \Phi(s):=\Phi(T-s).\]

Before we come to the results, we let us rewrite the stochastic integrals in a way that fits to our purposes. The integrals with respect to $dZ_s$ can be written as $\hat\pi$-integrals. To this end, let us define mappings
$E,\;F,\;G,\;\tilde F$ from $[0,T]\times\mathcal O\times\R$ into $H$, by
\begin{eqnarray}
E(s,\xi,\sigma)&:=& \sum_{k=1}^\infty e_k(\xi)\sigma S(T-s)Q^{1/2}e_k,\\
F(s,\xi,\sigma)&:=&\sum_{k=1}^\infty e_k(\xi)\sigma \Phi(s)e_k,\\
G(s,\xi,\sigma)&:=&\sum_{k=1}^\infty e_k(\xi)\sigma\Gamma(s)e_k,\label{G}\\
\tilde F(s,\xi,\sigma)&:=&\sum_{k=1}^\infty e_k(\xi)\sigma \tilde \Phi(s)e_k,
\end{eqnarray}
where $(e_k)_{k\in\N}$ is an orthonormal basis of $H$ and the infinite sums are limits in the space $L^2([0,T]\times\mathcal O\times\R,\,ds\,d\xi\,\nu(d\sigma);\,H)$. These limits exist since the operator valued functions $s\mapsto S(s)Q^{1/2},\;s\mapsto\Phi(s)$ and $s\mapsto\Gamma(s)$ belong to $L^2\big([0,T],\,ds;\,L_{\text{(HS)}}(H)\big)$ and because of the integrability assumption (\ref{AssNu0}). As usual, we do not distinguish between measurable mappings $[0,T]\times\mathcal O\times\R\to H$  and their $ds\,d\xi\,\nu(d\sigma)$-equivalence classes.
Notice that the following equalities hold:
\begin{eqnarray*}
\int_0^tS(T-s)Q^{1/2}\,dZ_s&=& \int_0^t\int_{\mathcal O\times\R}E(s,\xi,\sigma)\,\hat\pi(ds,d\xi,d\sigma),\\
\int_0^t\Phi(s)\,dZ_s &=&\int_0^t\int_{\mathcal O\times\R}F(s,\xi,\sigma)\,\hat\pi(ds,d\xi,d\sigma),\\
\int_0^t\Gamma(s)\,dZ_s &=&\int_0^t\int_{\mathcal O\times\R}G(s,\xi,\sigma)\,\hat\pi(ds,d\xi,d\sigma),\\
\int_0^t\tilde \Phi(s)\,dZ_s &=&\int_0^t\int_{\mathcal O\times\R}\tilde F(s,\xi,\sigma)\,\hat\pi(ds,d\xi,d\sigma).
\end{eqnarray*}
Indeed, if we consider for example $\int_0^t\Phi(s)\,dZ_s$, it is easy to verify that
\begin{eqnarray*}
\int_0^t\Phi(s)\,dZ_s&=& L^2(\Omega,\mathcal A,\wP;H)\text{-}\lim_{n\to\infty}\sum_{k=1}^n\int_0^t\Phi(s)e_k\,dZ^{(k)}_s\\
&=& L^2(\Omega,\mathcal A,\wP;H)\text{-}\lim_{n\to\infty}\sum_{k=1}^n \int_0^t\int_{\mathcal O\times\R}e_k(\xi)\sigma \Phi(s)e_k\,\hat\pi(ds,d\xi,d\sigma)\\
&=& \int_0^t\int_{\mathcal O\times\R}\left(\sum_{k=1}^\infty e_k(\xi)\sigma \Phi(s)e_k\right)\,\hat\pi(ds,d\xi,d\sigma),
\end{eqnarray*}
where the infinite sum in the last integral is a limit in $L^2([0,T]\times\mathcal O\times\R,\,ds\,d\xi\,\nu(d\sigma);\,H)$.

\section{Error expansion}

In this section, for suitable functions $\varphi$ defined on $H$, we state and prove representation formulas for the time discretization error $\E\varphi(X^N_h)-\E\varphi(X_{h,T})$ and the space discretization error $\E\varphi(X_{h,T})-\E\varphi(X_T)$. The errors are represented in terms of the functions $v_h:\,[0,T]\times H\to\R$ and $v_h:\,[0,T]\times H\to\R$ defined by
\begin{eqnarray}
v_h(t,x) &:=& \E\varphi\left( x+ \int_{T-t}^T S_h(T-r)P_hQ^{1/2}dZ_r\right),\label{defvh}\\
v(t,x) &:=& \E\varphi\left( x+ \int_{T-t}^T S(T-r)Q^{1/2}dZ_r\right).\label{defv}
\end{eqnarray}
For the reader's convenience, we summarize all equations, definitions and assumptions which will be needed in the formulation of our error representation theorem below:
\bigskip\\
\begin{samepage}
{\bf Summary/Assumptions:}
\begin{itemize}
\item The equation we are interested in is
\begin{equation}\tag{\ref{sde}}
dX_t + AX_t\,dt =  Q^{1/2}\, dZ_t,\quad X_0=x_0\in H,\quad t\in [0,T].
\end{equation}
\item \emph{Assumption:} The jump size intensity $\nu$ of $(Z_t)_{t\in[0,T]}$ satisfies
\begin{equation}\tag{\ref{AssNu0}}
\int_\R\sigma^2\,\nu(d\sigma)<\infty.
\end{equation}
\item \emph{Assumption:} The mild solution of equation \eqref{sde} exists, i.e. \begin{equation}\tag{\ref{intCond}}
    \|S(t)Q^{1/2}\|_{\text{(HS)}}\in L^2([0,T],dt)
    \end{equation}
    (This ensures that the weak solution exists and is equal to the mild solution.)
\item We consider the following space-time- and space-discretizations of the weak solution of \eqref{sde}
\begin{align}
X_h^n &= S_{h,\Delta t}^n P_hx_0 +\sum_{k=0}^{n-1}S_{h,\Delta t}^{n-k-1}T_{h,\Delta t}P_h Q^{1/2}(Z_{t_{k+1}}-Z_{t_k}),\quad n\in\{0,\ldots,N\},\tag{\ref{discrSoln}}\\
\tag{\ref{Xht}}
X_{h,t} &= S_h(t)P_hx_0 + \int_0^tS_h(t-s)P_hQ^{1/2}\,dZ_s,\qquad t\in [0,T].
\end{align}
\end{itemize}
\end{samepage}

\begin{thm}\label{errorRep}
Assume that \eqref{AssNu0} and \eqref{intCond} hold. Let $\varphi\in C^2_b(H)$ such that for each $x\in H$ the derivative $D^2\varphi(x)$ is an element of $L_{\text{(HS)}}(H)$ and that the mapping $x\mapsto D^2\varphi(x)\in L_{\text{(HS)}}(H)$ is uniformly continuous on any bounded subset of $H$. Let $T\geq 1$ and $(X_t)_{t\in[0,T]}$ be the $H$-valued stochastic process satisfying equation \eqref{sde}. For any $N\geq 1$ and $h>0$, let $(X^n_h)_{n\in\{0,\ldots, N\}}$ be given by $(\ref{discrSoln})$ and let $(X_{h,t})_{t\in[0,T]}$ be as in $(\ref{Xht})$.\\
Then the following error expansions hold:
\begin{align}
\E\varphi(X^N_h)-\E\varphi(X_{h,T})
&=\left\{v_h(T,S_{h,\Delta t}^NP_hx_0)-v_h(T,S_h(T)P_hx_0)\right\}\notag\\
&+\E\int_0^T\int_{\mathcal O\times \R}\Big\{v_h\big(T-t,\bar Y_{h,t-}+G(t,\xi,\sigma)\big)-v_h\big(T-t,\bar Y_{h,t-}+F(t,\xi,\sigma)\big)\notag\\
&\qquad\qquad\quad\;+\left\langle D_xv_h(T-t,\bar Y_{h,t-}),\,F(t,\xi,\sigma)-G(t,\xi,\sigma)\right\rangle_H\Big\}\,dt\,d\xi\,\nu(d\sigma)\notag\\
&=: I + II,\label{o}
\end{align}
\begin{align}
\E\varphi(X_{h,T})-\E\varphi(X_T)
&=\left\{v(T,S_h(T)P_hx_0)-v(T,S(T)x_0)\right\}\notag\\
&+\E\int_0^T\int_{\mathcal O\times \R}\Big\{v\big(T-t, Y_{h,t-}+F(t,\xi,\sigma)\big)-v\big(T-t,Y_{h,t-}+E(t,\xi,\sigma)\big)\notag\\
&\qquad\qquad\quad+\left\langle D_xv(T-t, Y_{h,t-}),\,E(t,\xi,\sigma)-F(t,\xi,\sigma)\right\rangle_H\Big\}\,dt\,d\xi\,\nu(d\sigma)\notag\\
&=: III + IV.\label{VI}
\end{align}
\end{thm}

\begin{rem}
In the following proof of Theorem \ref{errorRep}, the uniform continuity of the mapping $H\ni x\mapsto D^2\varphi(x)\in L_{\text{(HS)}}(H)$ on bounded subsets of $H$ is needed to be able to apply Itô's formula in infinite dimensions, see \cite{Met}, Section 27. This assumption is always fulfilled in finite dimensions and in literature it is sometimes forgotten to mention in the infinite-dimensional case.
\end{rem}

\begin{proof}[Proof of Theorem \ref{errorRep}]
We begin with the time discretization error $\E\varphi(X^N_h)-\E\varphi(X_{h,T})$. Due to the definition of $v_h$ we have
\begin{eqnarray}
\E\varphi(X_{h,T}) &=& v_h(T,S_h(T)P_hx_0),\label{x}\\
\E\varphi(X^N_h) &=& \E\varphi(\bar Y_{h,T}) \;\;=\;\;\E v_h(0,\bar Y_{h,T})\label{xx}.
\end{eqnarray}
Applying Itô's formula to the function $(t,x)\mapsto v_h(T-t,x)$ and the (càdlàg) martingale $(\bar Y_{h,t})_{t\in[0,T]}$ yields
\begin{equation}
\begin{aligned}\label{xxx}
v_h(0,\bar Y_{h,T}) &=
 v_h(T,\bar Y_{h,0}) - \int_0^T\frac {\partial v_h}{\partial t}(T-t,\bar Y_{h,t-})\,dt\\
&\quad+\int_0^T\left\langle D_xv_h(T-t,\bar Y_{h,t-}),\,d\bar Y_{h,t}\right\rangle_H
+\,\frac 12\int_0^T\left\langle D^2_xv_h(T-t,\bar Y_{h,t-}),\,d\llbracket\bar Y_{h,\cdot}\rrbracket_t^c\right\rangle_{\text{(HS)}}\\
&\quad+\sum_{t\leq T}\Big\{v_h(T-t,\bar Y_{h,t})-v_h(T-t,\bar Y_{h,t-}) -\left\langle D_xv_h(T-t,\bar Y_{h,t-}),\,\Delta\bar Y_{h,t}\right\rangle_H\Big\}.
\end{aligned}
\end{equation}
A process is called càdlàg (continu à droite et limites à gauche) if almost every path is right-continuous and has finite left limits. We have used the standard notation $\bar Y_{h,t-}=\lim_{s\nearrow t}\bar Y_{h,s}$ and $\Delta\bar Y_{h,t}=\bar Y_{h,t}-\bar Y_{h,t-}$ which make sense for every càdlàg process.
We will now consider the mean values of the terms on the right hand side of (\ref{xxx}) separately:
\begin{itemize}
\item[(i)] $\E v_h(T,\bar Y_{h,0})=\E v_h(T,S_{h,\Delta t}^NP_h x_0)=v_h(T,S_{h,\Delta t}^NP_hx_0)$.
\item[(ii)] $\E\int_0^T\left\langle D_xv_h(T-t,\bar Y_{h,t-}),\,d\bar Y_{h,t}\right\rangle_H=0$ because $\left(\int_0^t\left\langle D_xv_h(T-s,\bar Y_{h,s-}),\,d\bar Y_{h,s}\right\rangle_H\right)_{t\in[0,T]}$ is a martingale starting in $0$.
\item[(iii)] $\left(\llbracket\bar Y_{h,\cdot}\rrbracket^c_t\right)_{t\in[0,T]}=0$ because $(\bar Y_{h,t})_{t\in[0,T]}$ is a purely discontinuous martingale. This follows because $(\bar Y_{h,t})_{t\in[0,T]}$ is the $\mathcal M^2_T(H)$-limit of finite sums of its coordinate processes which can be easily identified as purely discontinuous martingales. Since the space of all purely discontinuous $L^2(\wP)$-martingales $\mathcal M^{2,d}_T(H)$ is closed in $\mathcal M^2_T(H)$, the process $(\bar Y_{h,t})_{t\in[0,T]}$ is also purely discontinuous (cf. \cite{Met}, Chapter 4). Hence $\left(\llbracket\bar Y_{h,\cdot}\rrbracket^c_t\right)_{t\in[0,T]}=0$, and therefore
    \[\frac 12\E\int_0^T\left\langle D^2_xv_h(T-t,\bar Y_{h,t-}),\,d\llbracket\bar Y_{h,\cdot}\rrbracket_t^c\right\rangle_{\text{(HS)}}=0.\]
\item[(iv)] Concerning the jump term
    \begin{equation}\label{jump-sum}
    \E\sum_{t\leq T}\Big\{v_h(T-t,\bar Y_{h,t})-v_h(T-t,\bar Y_{h,t-}) -\left\langle D_xv_h(T-t,\bar Y_{h,t-}),\,\Delta\bar Y_{h,t}\right\rangle_H\Big\}
    \end{equation}
    we give only a brief sketch how to deal with it; the full argument is postponed to the appendix.\\
    Let $\mathfrak m$ be the jump counting measure of $(\bar Y_{h,t})_{t\in[0,T]}$, i.e.
    \[\mathfrak m((0,t]\times B)(\omega)=\sum_{s\leq t}\one_B(\Delta\bar Y_{h,s}(\omega)),\;\qquad t\in(0,T],\;B\in\mathcal B(H),\;\omega\in\Omega.\]
    The basic idea is to write the sum in ($\ref{jump-sum}$) as a (pathwise) integral with respect to $\mathfrak m$ and then to consider the random measure $\mathfrak m$ as an image measure of $\pi$ with respect to the function $G$ defined by $(\ref{G})$:
    \begin{eqnarray}
    \lefteqn{\E\sum_{t\leq T}\Big\{v_h(T-t,\bar Y_{h,t})-v_h(T-t,\bar Y_{h,t-}) -\left\langle D_xv_h(T-t,\bar Y_{h,t-}),\,\Delta\bar Y_{h,t}\right\rangle_H\Big\}}\notag\\
    &=& \E\int_0^T\int_H\Big\{v_h(T-t,\bar Y_{h,t-}+y)-v_h(T-t,\bar Y_{h,t-})-\left\langle D_xv_h(T-t,\bar Y_{h,t-}),\,y\right\rangle_H\Big\}\,\mathfrak m(dt,dy)\notag\\
    &=& \E\int_0^T\int_{\mathcal O\times \R}\Big\{v_h\big(T-t,\bar Y_{h,t-}+G(t,\xi,\sigma)\big)-v_h(T-t,\bar Y_{h,t-})\notag\\
    &&\qquad\qquad\qquad\qquad\qquad-\left\langle D_xv_h(T-t,\bar Y_{h,t-}),\,G(t,\xi,\sigma)\right\rangle_H\Big\}\,\pi(dt,d\xi,d\sigma).\label{jump-sum-trans}
\end{eqnarray}
Here the last term is the expectation of an stochastic $L^1$-integral with respect to the random measure $\pi$, cf. \cite{PesZab}, Section 8.7. Since $dt\,d\xi\,\nu(d\sigma)$ is the compensator of $\pi$,
\begin{eqnarray*}
\lefteqn{\E\int_0^T\int_{\mathcal O\times \R}\Big\{v_h\big(T-t,\bar Y_{h,t-}+G(t,\xi,\sigma)\big)-v_h(T-t,\bar Y_{h,t-})}\notag\\
    &&\qquad\qquad\qquad\qquad\qquad-\left\langle D_xv_h(T-t,\bar Y_{h,t-}),\,G(t,\xi,\sigma)\right\rangle_H\Big\}\,\pi(dt,d\xi,d\sigma)\notag\\
    &=& \E\int_0^T\int_{\mathcal O\times \R}\Big\{v_h\big(T-t,\bar Y_{h,t-}+G(t,\xi,\sigma)\big)-v_h(T-t,\bar Y_{h,t-})\notag\\
    &&\qquad\qquad\qquad\qquad\qquad-\left\langle D_xv_h(T-t,\bar Y_{h,t-}),\,G(t,\xi,\sigma)\right\rangle_H\Big\}\,dt\,d\xi\,\nu(d\sigma).
    \end{eqnarray*}
\item[(v)] Going back to the usual construction of the stochastic integral, one can see that the laws of the random variables
\[\int_{T-t}^T S_h(T-s)P_hQ^{1/2}\,dZ_s \;\;\text{ and }\; \int_0^t S_h(s)P_hQ^{1/2}\,dZ_s\]
are equal. Consequently
\begin{eqnarray*}
v_h(t,x)&=&\E\varphi\left(x+\int_0^t S_h(s)P_hQ^{1/2}dZ_s\right)\\
&=& \E\varphi\left(x+\int_0^t\tilde \Phi(s)\,dZ_s\right),
\end{eqnarray*}
where we have used the notation of Section \ref{ApproximationScheme} $\tilde\Phi(s)=\Phi(T-s),\;\Phi(s)=S_h(T-s)P_hQ^{1/2}$ for $s\in[0,T].$
Now we apply the Itô formula to the function $H\ni y\mapsto \varphi(x+y)\in\R$ and the (càdlàg) martingale $\left(\int_0^t\tilde \Phi(r)\,dZ_r\right)_{t\in[0,T]}$:
\begin{eqnarray*}
\varphi\left(x+\int_0^t\tilde \Phi(s)\,dZ_s\right)&=& \varphi(x) + \int_0^t\left\langle D\varphi\left(x+\int_0^{r-}\tilde \Phi(s)\,dZ_s\right),\,\tilde \Phi(r)\,dZ_r\right\rangle_H\\
&&+ \frac 12 \int_0^t\left\langle D^2\varphi\left(x+\int_0^{r-}\tilde \Phi(s)\,dZ_s\right),\,\tilde \Phi(r)^{\otimes2}\,d\llbracket Z\rrbracket^c_r\right\rangle_{\text{(HS)}}\\
&&+ \sum_{r\leq t}\left\{\varphi\left(x+\int_0^r\tilde \Phi(s)\,dZ_s\right)-\varphi\left(x+\int_0^{r-}\tilde \Phi(s)\,dZ_s\right)\right.\\
&&\left.\qquad\qquad-\left\langle D\varphi\left(x+\int_0^{r-}\tilde \Phi(s)\,dZ_s\right),\,\Delta \int_0^r\tilde \Phi(s)\,dZ_s\right\rangle_H\right\}.
\end{eqnarray*}
Here, $\int_0^{r-}\tilde \Phi(s)\,dZ_s$ denotes the (pathwise) left limit $\lim_{q\nearrow r}\int_0^q\tilde\Phi(s)\,dZ_s$ and the expression $\tilde \Phi(r)^{\otimes 2}$ in the second line of the formula stands for the operator $L_{\text{(HS)}}(H)\to L_{\text{(HS)}}(H),\;T\mapsto \tilde \Phi(r)T\tilde \Phi(r)^*$.
Taking the expectation on both sides of the equation, reasoning as in (iv) and differentiating with respect to $t$ yields
\begin{eqnarray}
\frac{\partial}{\partial t}v_h(t,x)&=&\int_{\mathcal O\times\R}\left\{v_h\big(t,x+\tilde F(t,\xi,\sigma)\big)-v_h(t,x)\right.\notag\\
&&\qquad\qquad\qquad-\left.\left\langle D_xv_h(t,x),\,\tilde F(t,\xi,\sigma)\right\rangle_H\right\}\,d\xi\,\nu(d\sigma).
\end{eqnarray}
\end{itemize}
Finally we get (\ref{o}) if we combine (\ref x) - (\ref{xxx}) and (i) - (v).
\bigskip

Now consider the spatial error $\E\varphi(X_{h,T})-\E\varphi(X_T)$. The definition of the function $v$ implies
\begin{eqnarray*}
\E\varphi(X_T) &=& v(T,S(T)x_0),\\
\E\varphi(X_{h,T})&=&\E v(0,Y_{h,T}).
\end{eqnarray*}
As above, we can apply Itô's formula to $(t,x)\mapsto v(T-t,x)$ and $(Y_{h,t})_{t\in[0,T]}$ to get
\begin{eqnarray*}
v(0,Y_{h,T}) &=& v(T,Y_{h,0}) - \int_0^T\frac {\partial v}{\partial t}(T-t, Y_{h,t-})\,dt\label{1}\\
&&+\int_0^T\left\langle D_xv(T-t,Y_{h,t-}),\,dY_{h,t}\right\rangle_H\notag\\
&&+\frac 12\int_0^T\left\langle D^2_xv(T-t,Y_{h,t-}),\,d\llbracket Y_{h,\cdot}\rrbracket_t^c\right\rangle_{\text{(HS)}}\notag\\
&&+\sum_{t\leq T}\Big\{v(T-t,Y_{h,t})-v(T-t,Y_{h,t-}) -\left\langle D_xv(T-t, Y_{h,t-}),\,\Delta Y_{h,t}\right\rangle_H\Big\}\notag.
\end{eqnarray*}
This can be used, similar to the argumentation leading to (\ref{o}), to verify $(\ref{VI})$.
\end{proof}\bigskip

\begin{rem}
The proof of Theorem~\ref{errorRep} reveals that we can extend Theorem~\ref{errorRep} to equations of the form
\begin{equation}\label{mixed}
dX_t + AX_t\,dt =  Q^{1/2}_0\,dW_t + Q_1^{1/2}\, dZ_t,\quad X_0=x_0\in H,\quad t\in [0,T],
\end{equation}
where $(W_t)_{t\in[0,T]}$ is a cylindrical Wiener process on $H$ which is independent of $(Z_t)_{t\in[0,T]}$ and the covariance operators $Q_0$ and $Q_1$ are, just as $Q$ above, bounded, nonnegative definite, symmetric and satisfy \eqref{ass1}, \eqref{assQ}. The solution $(X_t)_{t\in[0,T]}$ of $(\ref{mixed})$ is given by
\[X_t =S(t)x_0 +\int_0^tS(t-s)Q_1^{1/2}\,dW_s+\int_0^tS(t-s)Q_2^{1/2}\,dZ_s,\qquad t\in[0,T],\]
and the discretizations are
\begin{eqnarray*}
X_h^n &=& S_{h,\Delta t}^n P_hx_0 +\sum_{k=0}^{n-1}S_{h,\Delta t}^{n-k-1}T_{h,\Delta t}P_h Q_0^{1/2}(W_{t_{k+1}}-W_{t_k})\\
&&+\sum_{k=0}^{n-1}S_{h,\Delta t}^{n-k-1}T_{h,\Delta t}P_h Q_1^{1/2}(Z_{t_{k+1}}-Z_{t_k}),\qquad h>0,\; n\in\{0,\ldots,N\},\\
X_{h,t} &=& S_h(t)P_hx_0 + \int_0^tS_h(t-s)P_hQ_0^{1/2}\,dW_s\\
&&+\int_0^tS_h(t-s)P_hQ_1^{1/2}\,dZ_s,\qquad t\in [0,T].
\end{eqnarray*}
Assuming for simplicity that  $Q_1=Q_2$, the time discretization error $\E\varphi(X^N_h)-\E\varphi(X_{h,T})$ and the spatial error $\E\varphi(X_{h,T})-\E\varphi(X_T)$ are obtained by adding the terms
\begin{equation}\label{WII}
\frac 12 \,\E\int_0^T\text{Tr}\Big\{\Gamma(t)^*\,D^2_xv_h(T-t,\bar Y_{h,t-})\,\Gamma(t)\,-\,
\Phi(t)^*D^2_xv_h(T-t,\bar Y_{h,t-})\Phi(t)\Big\}\,dt
\end{equation}
and
\begin{equation}\label{WIV}
\frac 12 \,\E\int_0^T\text{Tr}\Big\{\Phi(t)^*\,D^2_xv(T-t, Y_{h,t-})\,\Phi(t)\,-\,
\big(S(T-t)Q_1^{1/2}\big)^*D^2_xv(T-t,Y_{h,t-})\big(S(T-t)Q_1^{1/2}\big)\Big\}\,dt
\end{equation}
to the right hand side of $(\ref{o})$ and $(\ref{VI})$, respectively. (Here one has to replace $Q$ by $Q_1$ in the definition of $\Phi(t)$ and $\Gamma(t)$.)
\end{rem}

\section{Weak order of convergence}
In this section, we show an estimate for the weak order of the convergence of $X_h^N$ to $X_T$ as $h$ tends to zero and $N$ tends to infinity, given the integrability condition (\ref{AssNu}). The proof is based on the error expansions in the last section.\\
For the convenience of the readers we summarize the assumptions made in the theorem below.\bigskip\\
\begin{samepage}
{\bf Assumptions:}
\begin{itemize}
\item The jump size intensity $\nu$ of $(Z_t)_{t\in[0,T]}$ satisfies
\begin{equation}\tag{\ref{AssNu}}
\int_\R\max\big(|\sigma|,\sigma^2\big)\,\nu(d\sigma)<\infty.
\end{equation}
\item The finite element spaces $V_h,\;h>0,$ are such that for all $q\in[0,2]$ there exist constants $\kappa_1>0,\;\kappa_2>0$ independent from $h$ such that for all $t>0$
\begin{align}
\|S_h(t)P_h-S(t)\|_{L(H)}&\leq \kappa_1 h^qt^{-q/2},\tag{\ref{S1}}\\
\|S_h(t)P_h-S(t)\|_{L(H,D(A^{1/2}))}&\leq \kappa_2 ht^{-1}.\tag{\ref{S2}}
\end{align}
\item There exist real numbers $\alpha$ and $\beta$ with
\begin{gather}
\alpha>0,\quad\beta\in(\alpha-1,\alpha],\tag{\ref{ass1}}\\
\text{Tr}(A^{-\alpha})=\sum_{n=1}^\infty\lambda_n^{-\alpha}< \infty\;,\tag{\ref{assTr}}\\
A^\beta Q\in L(H).\tag{\ref{assQ}}
\end{gather}
(Note that this particularly implies assumption \eqref{intCond} needed for Theorem~\ref{errorRep}.)
\end{itemize}
\end{samepage}
\begin{thm}\label{result}
Assume \eqref{AssNu}, \eqref{S1}, \eqref{S2}, \eqref{ass1}, \eqref{assTr} and $\eqref{assQ}$ listed above. Let $\varphi\in C^2_b(H)$ such that for each $x\in H$ the derivative $D^2\varphi(x)$ is an element of $L_{\text{(HS)}}(H)$ and that the mapping $x\mapsto D^2\varphi(x)\in L_{\text{(HS)}}(H)$ is uniformly continuous on any bounded subset of $H$. Let $T\geq 1$ and $(X_t)_{t\in[0,T]}$ be the $H$-valued stochastic process satisfying the equation
\begin{equation}\tag{\ref{sde}}
dX_t + AX_t\,dt =  Q^{1/2}\, dZ_t,\quad X_0=x_0\in H,\quad t\in [0,T].
\end{equation}
For any $N\geq 1$ and $h>0$, let $(X^n_h)_{n\in\{0,\ldots, N\}}$ be given by $(\ref{discrSoln})$.\\
Then for any $\gamma<1-\alpha+\beta\leq 1$, there exists a constant $C=C(T,\varphi,A,Q,|x_0|,\gamma,\nu)>0$ which does not depend on $h$ and $N$ such that the following inequality holds
\[ |\E\varphi(X^N_h)-\E\varphi(X_T)|\leq C\cdot (h^{2\gamma}+(\Delta t)^{\gamma}),\]
where $\Delta t =T/N\leq 1$.
\end{thm}
\begin{proof}
As in Theorem~\ref{errorRep}, we split the error into the time discretization error and the spatial error,
\[\E\varphi(X^N_h)-\E\varphi(X_T)=\Big\{\E\varphi(X^N_h)-\E\varphi(X_{h,T})\Big\}+
\Big\{\E\varphi(X_{h,T})-\E\varphi(X_T)\Big\}.\]
We will estimate each term separately. Throughout this proof, $C$ denotes a positive (and finite) constant that may change from line to line.

Consider the time discretization error given by $(\ref{o})$, $\E\varphi(X^N_h)-\E\varphi(X_{h,T})=I+II$.
Clearly, due to the definition of $v_h$, \eqref{defvh},
\[|I| \leq \|\varphi\|_{C^1_b(H)}\left\|S^N_{h,\Delta t}P_h-S_h(T)P_h\right\|_{L(H)}|x_0|_H.\]
Using spectral calculus and $(\ref{theta})$, $\left\|S^N_{h,\Delta t}P_h-S_h(T)P_h\right\|_{L(H)}$ can be bounded uniformly with respect to $h$. To be more precise, for $T\geq 1$ we have
\begin{eqnarray*}
\left\|S^N_{h,\Delta t}P_h-S_h(T)P_h\right\|_{L(H)}&\leq&\sup_{\lambda>0}\left|e^{-N\lambda\Delta t}-\left(\frac{1-(1-\theta)\lambda\Delta t}
{1+\theta\lambda\Delta t}\right)^N\right|\\
&=&\sup_{r>0}\left|e^{-Nr}-\left(\frac{1-(1-\theta)r}
{1+\theta r}\right)^N\right|\\
&\leq& \frac CN\;\;\leq\;\; C\cdot\Delta t,
\end{eqnarray*}
cf. \cite{LeRoux}, p. 921, Theorem 1.1 for the penultimate estimate.
One obtains
\begin{equation}\label{Roux}
|I|\leq C\cdot\Delta t.
\end{equation}
Concerning $II$, we apply the mean value theorem and the Cauchy-Schwarz inequality and obtain
\begin{eqnarray}
II&\leq& \int_0^T\int_{\mathcal O\times\R}2\|D_xv_h\|_{C_b(H)}\left|F(t,\xi,\sigma)-G(t,\xi,\sigma)\right|\,dt\,d\xi\,\nu (d\sigma)\notag\\
&=&2\|D_xv_h\|_{C_b(H)}\int_\R|\sigma|\,\nu(d\sigma)\int_0^T\int_{\mathcal O}\left|\sum_{k=1}^\infty e_k(\xi)\big(\Phi(t)-\Gamma(t)\big)e_k\right|_H\,dt\,d\xi\label{*},
\end{eqnarray}
where the expression $\sum_{k=1}^\infty e_k(\xi)\big(\Phi(t)-\Gamma(t)\big)e_k$ has to be understood as the value at $(t,\xi)$ of a $dt\,d\xi$-version of \[L^2([0,T]\times\mathcal O, dt\,d\xi;\,H)\text{-}\lim_{N\to\infty}\left((t,\xi)\mapsto \sum_{k=1}^N e_k(\xi)\big(\Phi(t)-\Gamma(t)\big)e_k\right).\]
Next, let $\gamma>0$ and $\gamma_1>0$ such that $0<\gamma<\gamma_1<1-\alpha+\beta\leq1$. Using
\begin{eqnarray*}
\lefteqn{\sum_{k=1}^\infty e_k(\xi)\big(\Phi(t)-\Gamma(t)\big)e_k=}\\
&&\left(S_h(T-t)-\sum_{n=0}^{N-1}\one_{(t_n,t_{n+1}]}(t)S_{h,\Delta t}^{N-n-1}T_{h,\Delta t}\right)A_h^{(1-\gamma_1)/2}\sum_{k=1}^\infty e_k(\xi)A_h^{(\gamma_1-1)/2}P_hQ^{1/2}e_k,
\end{eqnarray*}
(\ref{*}) can be estimated by
\begin{align}
II&\leq C\int_0^T\left\|\left(S_h(T-t)-\sum_{n=0}^{N-1}\one_{(t_n,t_{n+1}]}(t)S_{h,\Delta t}^{N-n-1}T_{h,\Delta t}\right)A_h^{(1-\gamma_1)/2}\right\|_{L(V_h)}\,dt\notag\\
&\qquad\qquad\qquad\times\int_{\mathcal O}\left|\sum_{k=1}^\infty e_k(\xi)A_h^{(\gamma_1-1)/2}P_hQ^{1/2}e_k\right|_H\,d\xi\notag\\
&\leq C\int_0^T\left\|\left(S_h(T-t)-\sum_{n=0}^{N-1}\one_{(t_n,t_{n+1}]}(t)S_{h,\Delta t}^{N-n-1}T_{h,\Delta t}\right)A_h^{(1-\gamma_1)/2}\right\|_{L(V_h)}\,dt\notag\\
&\qquad\qquad\qquad\times\left(\int_{\mathcal O}\one_{\mathcal O}(\xi)\,d\xi\right)^{1/2}\left(\int_{\mathcal O}\left|\sum_{k=1}^\infty e_k(\xi)A_h^{(\gamma_1-1)/2}P_hQ^{1/2}e_k\right|_H^2\,d\xi\right)^{1/2}\notag\\
&= C\,\|A_h^{(\gamma_1-1)/2}P_hQ^{1/2}\|_{\text{(HS)}}\notag\\
&\qquad\qquad\times\int_0^T\left\|\left(S_h(T-t)-\sum_{n=0}^{N-1}\one_{(t_n,t_{n+1}]}(t)S_{h,\Delta t}^{N-n-1}T_{h,\Delta t}\right)A_h^{(1-\gamma_1)/2}\right\|_{L(V_h)}\,dt.\label{**}
\end{align}
Following \cite{DebPrin}, $\big\|A_h^{(\gamma_1-1)/2}P_hQ^{1/2}\big\|_{\text{(HS)}}$ can be bounded from above by some constant $C>0$ which does not depend on the choice of $h>0$:
\begin{equation}\label{i}
\big\|A_h^{(\gamma_1-1)/2}P_hQ^{1/2}\big\|_{\text{(HS)}}\leq C.
\end{equation}
This is essentially due to (\ref{assTr}) and (\ref{assQ}). Furthermore, again following \cite{DebPrin},
\begin{eqnarray}
\lefteqn{\left\|\left(S_h(T-t)-\sum_{n=0}^{N-1}\one_{(t_n,t_{n+1}]}(t)S_{h,\Delta t}^{N-n-1}T_{h,\Delta t}\right)A_h^{(1-\gamma_1)/2}\right\|_{L(V_h)}}\notag\\
&&\leq\quad
\begin{cases}
\quad C\Delta t^\gamma\sum_{n=0}^{N-2}\one_{(t_n,t_{n+1}]}(t)((N-n-1)\Delta t)^{-((1-\gamma_1)/2+\gamma)}&\quad t\in(0,t_{N-1}]\\
\quad C(T-t)^{(\gamma_1-1)/2},&\quad t\in(t_{N-1},T],
\end{cases}\label{ii}
\end{eqnarray}
which follows directly from spectral calculus.
A combination of (\ref{**}), (\ref{i}) and (\ref{ii}) yields
\begin{equation}\label{iii}
\begin{aligned}
II&\leq C\cdot\left(\Delta t^\gamma\int_0^{t_{N-1}}\sum_{n=0}^{N-2}\one_{(t_n,t_{n+1}]}(t)((N-n-1)\Delta t)^{-((1-\gamma_1)/2+\gamma)}\,dt+\int_{t_{N-1}}^T(T-t)^{(\gamma_1-1)/2}\,dt\right)\\
&\leq C\cdot\Delta t^\gamma,
\end{aligned}
\end{equation}
as $(1-\gamma_1)/2+\gamma\in(0,1)$ and $(T-t)^{(\gamma_1-1)/2}\leq(T-t)^{-((1-\gamma_1)/2+\gamma)}\Delta t^\gamma$ for $t\in(t_{N-1},T]$.

Finally (\ref{o}), (\ref{Roux}) and (\ref{iii}) add up to
\begin{equation}
\big|\E\varphi(X^N_h)-\E\varphi(X_{h,T})\big|\leq C\Delta t^\gamma
\end{equation}
for all $T\geq 1$ and $\Delta t\leq 1$.
\bigskip

Now we turn to the spatial error given by $(\ref{VI})$, $\E\varphi(X_{h,T})-\E\varphi(X_T)=III+IV$.
Firstly, according to the definition of $v$, \eqref{defv}, and resulting from (\ref{S1}) with $q=2\gamma<2$,
\begin{align}\label{V}
|III|&\;\leq\; \|\varphi\|_{C^2_b(H)}\left\|S_h(T)P_h-S(T)\right\|_{L(H)}|x_0|_H\notag\\
&\;\leq\; \|\varphi\|_{C^2_b(H)}\left(\kappa_1h^{2\gamma}T^{-\gamma}\right)|x_0|_H\;= \;C\cdot h^{2\gamma}.
\end{align}
Considering $IV$, we apply the mean value theorem and the Cauchy-Schwarz inequality and obtain \begin{eqnarray}
|IV|&\leq& \int_0^T\int_{\mathcal O\times\R}2\|D_xv\|_{C_b(H)}\left|E(t,\xi,\sigma)-F(t,\xi,\sigma)\right|\,dt\,d\xi\,\nu (d\sigma)\notag\\
&=&2\|D_xv\|_{C_b(H)}\int_\R|\sigma|\,\nu(d\sigma)\notag\\
&&\times\int_0^T\int_{\mathcal O}\left|\sum_{k=1}^\infty e_k(\xi)\big(S(T-t)-S_h(T-t)P_h\big)Q^{1/2}e_k\right|_H\,dt\,d\xi\label{2},
\end{eqnarray}
where the expression $\sum_{k=1}^\infty e_k(\xi)\big(S(T-t)-S_h(T-t)P_h\big)Q^{1/2}e_k$ has to be understood as the value at $(t,\xi)$ of a $dt\,d\xi$-version of \[L^2([0,T]\times\mathcal O, dt\,d\xi;\,H)\text{-}\lim_{N\to\infty}\left((t,\xi)\mapsto \sum_{k=1}^N e_k(\xi)\big(S(T-t)-S_h(T-t)P_h\big)Q^{1/2}e_k\right).\]
Pick $\gamma_1>0$ such that $0<\gamma<\gamma_1<1-\alpha+\beta\leq 1$. Because of ($\ref{interpol}$) we know that $A^{\beta/2}Q^{1/2}\in L(H)$. Furthermore, the fact that $1-\gamma_1+\beta>\alpha$ implies $A^{-(1-\gamma_1+\beta)/2}\in L_{\text{(HS)}}(H)$. It is also not hard to see that the operator $(S(T-t)-S_h(T-t)P_h)A^{(1-\gamma_1)/2}$ has a continuous extension defined on the whole space $H$. Therefore we may write
\begin{eqnarray*}
\lefteqn{\sum_{k=1}^\infty e_k(\xi)\big(S(T-t)-S_h(T-t)P_h\big)Q^{1/2}e_k}\\
&=& \left(S(T-t)-S_h(T-t)P_h\right)A^{(1-\gamma_1)/2} \sum_{k=1}^\infty e_k(\xi)A^{-(1-\gamma_1+\beta)/2}A^{\beta/2}Q^{1/2}e_k.
\end{eqnarray*}
A further application of the Cauchy-Schwarz inequality in (\ref{2}) gives
\begin{eqnarray}
|IV|&\leq&C\int_0^T\Big\|\big(S(T-t)-S_h(T-t)P_h\big)A^{(1-\gamma_1)/2}\Big\|_{L(H)}\,dt\notag\\
&&\qquad\qquad\qquad\times\left(\int_{\mathcal O}\one_{\mathcal O}(\xi)\,d\xi\right)^{1/2}\left(\int_{\mathcal O}\left|\sum_{k=1}^\infty e_k(\xi)A^{-(1-\gamma_1+\beta)/2}A^{\beta/2}Q^{1/2}e_k\right|_H^2\,d\xi\right)^{1/2}\notag\\
&=& C\,\big\|A^{-(1-\gamma_1+\beta)/2}A^{\beta/2}Q^{1/2}\big\|_{\text{(HS)}}
\int_0^T\Big\|\big(S(T-t)-S_h(T-t)P_h\big)A^{(1-\gamma_1)/2}\Big\|_{L(H)}\,dt\notag\\
&=& C\,\int_0^T\Big\|\big(S(T-t)-S_h(T-t)P_h\big)A^{(1-\gamma_1)/2}\Big\|_{L(H)}\,dt.\label{IV}
\end{eqnarray}

We are now going to show that for all $t\in[0,T]$
\begin{equation}\label{III}
\Big\|\left(S(T-t)-S_h(T-t)P_h\big)A^{(1-\gamma_1)/2}\right\|_{L(H)}\leq
C\,h^{2\gamma}\left(t^{-\big(\gamma_1(\gamma_1-1)/(2\gamma)+1\big)}+ t^{-\big((1-\gamma_1)/2+\gamma\big)}\right).
\end{equation}
This will be done in several steps.
\begin{itemize}
\item[(i)] Due to the self adjointness of $S(T-t),\;S_h(T-t)P_h$ and $A^{(1-\gamma_1)/2}$,
    \begin{eqnarray*}
    \Big\|\big(S(T-t)-S_h(T-t)P_h\big)A^{(1-\gamma_1)/2}\Big\|_{L(H)}&=&
    \Big\|\Big(\big(S(T-t)-S_h(T-t)P_h\big)A^{(1-\gamma_1)/2}\Big)^*\Big\|_{L(H)}\\
    &=&\Big\|A^{(1-\gamma_1)/2}\big(S(T-t)-S_h(T-t)P_h\big)\Big\|_{L(H)}.
    \end{eqnarray*}
    Here the operator $A^{(1-\gamma_1)/2}\big(S(T-t)-S_h(T-t)P_h\big)$ is properly defined since $S(t)H\subset D(A^\lambda),\;t>0,\;\lambda\in\R$ and $V_h\subset D(A^{1/2})\subset D(A^{(1-\gamma_1)/2})$.
\item[(ii)]
    Due to the Hölder inequality we have for all $x\in D(A^{1/2})$
    \[\big|A^{(1-\gamma_1)/2}x\big|_H\leq\big|A^{1/2}x\big|_H^{1-\gamma_1}\cdot|x|_H^{\gamma_1}.\]
    Consequently,
    \begin{eqnarray}
    \lefteqn{\big\|A^{(1-\gamma_1)/2}\big(S(T-t)-S_h(T-t)P_h\big)\big\|_{L(H)}}\notag\\
    &\leq&\|S(T-t)-S_h(T-t)P_h\|_{L(H,D(A^{1/2}))}^{1-\gamma_1}\cdot\|S(T-t)-S_h(T-t)P_h\|_{L(H)}^{\gamma_1}.\label{O}
    \end{eqnarray}
\item[(iii)]
    If $2\gamma+\gamma_1\geq 1$, then $(2\gamma+\gamma_1-1)/\gamma\in[0,2]$. In this case, we combine (\ref{O}) with (\ref{S2}) and (\ref{S1}), $q=(2\gamma+\gamma_1-1)/\gamma$. Provided that $h\leq 1$, one gets
    \begin{eqnarray*}
    \lefteqn{\big\|A^{(1-\gamma_1)/2}\left(S(T-t)-S_h(T-t)P_h\big)\right\|_{L(H)}}\\
    &\leq& C\,\big(h^{1-\gamma_1}t^{\gamma_1-1}\big)\big(h^{(2\gamma+\gamma_1-1)\gamma_1/\gamma}
    t^{-(2\gamma+\gamma_1-1)\gamma_1/(2\gamma)}\big)\\
    &\leq& C h^{2\gamma}t^{-\big(\gamma_1(\gamma_1-1)/(2\gamma)+1\big)}.
    \end{eqnarray*}
\item[(iv)] If $2\gamma+\gamma_1<1$, we first combine (\ref{O}) with (\ref{S2}) and (\ref{S1}) choosing $q=0$ and get
    \begin{equation}\label{I}
    \|A^{(1-\gamma_1)/2}\left(S(T-t)-S_h(T-t)P_h\big)\right\|_{L(H)}\leq
    C\,h^{1-\gamma_1}(T-t)^{-(1-\gamma_1)}.
    \end{equation}
    Secondly, again using (\ref{S1}) with $q=0$, one derives
    \begin{eqnarray}
    \|\lefteqn{A^{(1-\gamma_1)/2}\left(S(T-t)-S_h(T-t)P_h\big)\right\|_{L(H)}}\notag\\
    &\leq& \|S(T-t)-S_h(T-t)P_h\|_{L(H,D(A^{1/2}))}^{1-\gamma_1}
    \|S(T-t)-S_h(T-t)P_h\|_{L(H)}^{\gamma_1}\notag\\
    &\leq& C\, \|S(T-t)-S_h(T-t)P_h\|_{L(H,D(A^{1/2}))}^{1-\gamma_1}\notag\\
    &\leq& C\,\big( \|A^{1/2}S(T-t)\|_{L(H)}+\|A^{1/2}S_h(T-t)P_h\|_{L(H)}\big)^{1-\gamma_1}\notag\\
    &=& C\,\big( \|A^{1/2}S(T-t)\|_{L(H)}+\|A^{1/2}_hS_h(T-t)P_h\|_{L(H)}\big)^{1-\gamma_1}\notag\\
    &\leq& C\,(T-t)^{-(1-\gamma_1)/2}.\label{II}
    \end{eqnarray}
    In the last step, the inequality
    \[\sup_{x\geq0}x^\epsilon e^{-tx}\leq \left(\frac \epsilon e\right)^\epsilon t^{-\epsilon},\qquad t>0,\;\epsilon>0\]
    has been used.\\
    Now let $\lambda:=2\gamma/(1-\gamma_1)\in(0,1)$. Interpolating (\ref{I}) and (\ref{II})  yields
    \begin{eqnarray*}
     \lefteqn{\|A^{(1-\gamma_1)/2}\left(S(T-t)-S_h(T-t)P_h\big)\right\|_{L(H)}}\\
     &\leq&C\,\left(h^{1-\gamma_1}(T-t)^{-(1-\gamma_1)}\right)^\lambda
     \left((T-t)^{-(1-\gamma_1)/2}\right)^{1-\lambda}\\
     &=& C\,h^{2\gamma}(T-t)^{-((1-\gamma_1)/2+\gamma)}
    \end{eqnarray*}
\end{itemize}
Combining (i), (ii) and (iii) yields (\ref{III}).

As $\gamma_1(\gamma_1-1)/(2\gamma)+1<1$ and $(1-\gamma_1)/2+\gamma<1$, we can integrate (\ref{III}) with respect to the time variable $t$. Then (\ref{VI}), (\ref{V}), (\ref{IV}) and (\ref{III}) add up to
\[|\E\varphi(X_{h,T})-\E\varphi(X_T)|\leq C\,h^{2\gamma}.\]
\end{proof}

\begin{rem}
\begin{itemize}
\item[(i)]
The proof of Theorem~\ref{result} can easily be combined with the proof in \cite{DebPrin} in order to obtain the same result for the equation
\begin{equation}\label{mixed2}
dX_t + AX_t\,dt =  Q^{1/2}_0\,dW_t + Q_1^{1/2}\, dZ_t,\quad X_0=x_0\in H,\quad t\in [0,T],
\end{equation}
where $(W_t)_{t\in[0,T]}$ is a cylindrical Wiener process on $H$ which is independent of $(Z_t)_{t\in[0,T]}$ and where the covariance operators $Q_0$ and $Q_1$ are bounded, nonnegative definite, symmetric and satisfy \eqref{ass1}, \eqref{assQ}. The corresponding discretization and error expansion has been described in Remark 1.\\
Of course, the result also holds for the equation
\begin{equation}
dX_t + AX_t\,dt =  Q_0^{1/2}dW_t+Q^{1/2}_1\,dZ_{1,t} + Q_2^{1/2}\, dZ_{2,t},\quad X_0=x_0\in H,\quad t\in [0,T],
\end{equation}
where $(Z_{1,t})_{t\in[0,T]}$ and $(Z_{2,t})_{t\in[0,T]}$ are impulsive cylindrical processes satisfying the condition  required of $(Z_t)_{t\in[0,T]}$ (see \eqref{AssNu} above), the processes $(W_t)_{t\in[0,T]}$, $(Z_{1,t})_{t\in[0,T]}$, $(Z_{2,t})_{t\in[0,T]}$ are independent, and the covariance operators $Q_0$, $Q_1$ and $Q_2$ are bounded, nonnegative definite, symmetric and satisfy \eqref{ass1}, \eqref{assQ}. For example, one could consider impulsive cylindrical processes described by the jump size intensities
\[\nu_i(d\sigma)=\frac{1}{|\sigma|^{1+\alpha_i}}\one_{[-\tau,\tau]}(\sigma)\,d\sigma,\quad i=1,\,2,\]
with indices of stability $0<\alpha_1<\alpha_2<1$.
\item[(ii)]
One might try to avoid the integrability assumption (\ref{AssNu}) in the proof of Theorem~\ref{result} by rewriting the terms $II$ and $IV$ in (\ref{o}) and (\ref{VI}) using Taylor's theorem. Obviously,
\begin{eqnarray*}
II&=&\E\int_0^T\int_{\mathcal O\times\R}\Big\{\,\big\langle G(t,\xi,\sigma)\,,\,D^2_xv_h\big(T-t,\bar Y_{h,t-}+\theta G(t,\xi,\sigma)\big)G(t,\xi,\sigma)\big\rangle_H\\
&&\qquad\qquad -\,\big\langle F(t,\xi,\sigma)\,,\,D^2_xv_h\big(T-t,\bar Y_{h,t-}+\theta F(t,\xi,\sigma)\big)F(t,\xi,\sigma)\big\rangle_H\Big\}\,dt\,d\xi\,\nu(d\sigma),\\
IV&=&\E\int_0^T\int_{\mathcal O\times\R}\Big\{\,\big\langle F(t,\xi,\sigma)\,,\,D^2_xv\big(T-t,Y_{h,t-}+\vartheta F(t,\xi,\sigma)\big)F(t,\xi,\sigma)\big\rangle_H\\
&&\qquad\qquad -\,\big\langle E(t,\xi,\sigma)\,,\,D^2_xv\big(T-t, Y_{h,t-}+\vartheta E(t,\xi,\sigma)\big)E(t,\xi,\sigma)\big\rangle_H\Big\}\,dt\,d\xi\,\nu(d\sigma),
\end{eqnarray*}
where $\theta=\theta(t,\xi,\omega)\in (0,1)$ and $\vartheta=\vartheta(t,\xi,\omega)\in (0,1)$. However, the integrands appearing here cannot be estimated analogously to the estimates of the integrands in the terms (\ref{WII}) and (\ref{WIV}) in \cite{DebPrin}, which appear in the case of Gaussian noise. The reason is, that in the case of impulsive noise one has to estimate the difference of the second derivatives. As one needs integrability of suitable upper bounds for the integrands, this leads to profound technical complications.
\end{itemize}
\end{rem}

\begin{appendix}
\section{Proof of equality $(\ref{jump-sum-trans})$.}
Here we give a detailed proof for equality $(\ref{jump-sum-trans})$. We have to show
\begin{equation}\tag{\ref{jump-sum-trans}}
\begin{aligned}
    \lefteqn{\E\sum_{t\leq T}\Big\{v_h(T-t,\bar Y_{h,t})-v_h(T-t,\bar Y_{h,t-}) -\left\langle D_xv_h(T-t,\bar Y_{h,t-}),\,\Delta\bar Y_{h,t}\right\rangle_H\Big\}}\notag\\
   &= \E\int_0^T\int_{\mathcal O\times \R}\Big\{v_h\big(T-t,\bar Y_{h,t-}+G(t,\xi,\sigma)\big)-v_h(T-t,\bar Y_{h,t-})\notag\\
    &\qquad\qquad\qquad\qquad\qquad-\left\langle D_xv_h(T-t,\bar Y_{h,t-}),\,G(t,\xi,\sigma)\right\rangle_H\Big\}\,\pi(dt,d\xi,d\sigma).
\end{aligned}
\end{equation}
We have already seen that
\[
\bar Y_{h,t}-S^N_{h,\Delta t}P_hx_0=\int_0^t\Gamma(s)dZ_s
=\int_0^t\int_{\mathcal O\times \R} G(s,\xi,\sigma)\,\hat\pi(ds,d\xi,d\sigma),
\]
with $G(s,\xi,\sigma)=\sum_{k=1}^\infty e_k(\xi)\sigma\Gamma(s)e_k$, where the the infinite sum is a limit in the space $L^2([0,T]\times\mathcal O\times\R,\,ds\,d\xi\,\nu(d\sigma);\,H)$. By a standard monotone class argument, there exist simple functions
\[G_n(s,\xi,\sigma)=\sum_{k=1}^{m(n)}a_{n,k}\one_{A_{n,k}}(s,\xi,\sigma),\qquad (s,\xi,\sigma)\in\,[0,T]\times\mathcal O\times\R,\]
where $n\in\N,\;m(n)\in\N,\;a_{n,k}\in H$ and $A_{n,k}\in\mathcal B([0,T]\times\mathcal O\times\R)$ for $k=1,\ldots,m(n)$ with
\[\int_0^T\int_{\mathcal O\times\R}\one_{A_{n,k}}(t,\xi,\sigma)\,dt\,d\xi\,\nu(d\sigma)<\infty,\]
such that $G_n$ converges to $G$ in $L^2([0,T]\times\mathcal O\times\R,\,ds\,d\xi\,\nu(d\sigma);\,H)$ if $n$ tends to infinity. Note that the processes
\[\left(\int_0^t\int_{\mathcal O\times\R}G(s,\xi,\sigma)\,\hat\pi(ds,d\xi,d\sigma)\right)_{t\in[0,T]}\;\text{ and }\;\left(\int_0^t\int_{\mathcal O\times\R}G_n(s,\xi,\sigma)\,\hat\pi(ds,d\xi,d\sigma)\right)_{t\in[0,T]},\]
$n\in\N$, are martingales and hence have càdlàg modifications. Applying Doob's inequality for submartingales, we have
\begin{align}
\lefteqn{
\E\sup_{t\leq T}\left|\int_0^t\int_{\mathcal O\times\R}G(s,\xi,\sigma)\,\hat\pi(ds,d\xi,d\sigma)-
\int_0^t\int_{\mathcal O\times\R}G_n(s,\xi,\sigma)\,\hat\pi(ds,d\xi,d\sigma)\right|_H^2}\notag\\
&&\qquad\qquad\leq 4\,\E\left|\int_0^T\int_{\mathcal O\times\R}\left(G(s,\xi,\sigma)-G_n(s,\xi,\sigma)\right)\,\hat\pi(ds,d\xi,d\sigma)\right|_H^2\notag\\
&&\qquad\qquad=\;4\,\left|G-G_n\right|_{L^2([0,T]\times\mathcal O\times\R,\,ds\,d\xi\,\nu(d\sigma);\,H)}^2\label{glmKgz}\tag{A.1}
\end{align}
and the last term goes to zero as $n$ tends to infinity.\\
Remember that the Poisson random measure $\pi$ on $[0,T]\times\mathcal O\times\R$ with given compensator $dt\,d\xi\,\nu(d\sigma)$ is given by \[\pi=\sum_{j=1}^\infty\delta_{(T_j,\Xi_j,\Sigma_j)},\]
where $\big((T_j,\Xi_j,\Sigma_j)\big)_{j\in\N}$ is a properly chosen sequence of random elements in $[0,T]\times\mathcal O\times\R$, cf. \cite{PesZab}, Chapter 6.
Since the functions $G_n,\;n\in\N,$ are simple functions, the stochastic integral of $G_n$ with respect to $\hat\pi$ is just an $\omega$-wise integral with respect to
\[\sum_{j=1}^\infty\delta_{(T_j(\omega),\Xi_j(\omega),\Sigma_j(\omega))}(dt,d\xi,d\sigma)-dt\,d\xi\,\nu(d\sigma).\]
Therefore, for every $n\in\N$, the càdlàg process
\[\left(\int_0^t\int_{\mathcal O\times\R}G_n(s,\xi,\sigma)\,\hat\pi(ds,d\xi,d\sigma)\right)_{t\in[0,T]}\]
has the jumps $G_n\big(T_j(\omega),\Xi_j(\omega),\Sigma_j(\omega)\big)$, $j\in\N,\;T_j(\omega)\leq T,$ occurring at the jump times $T_j(\omega),\;j\in\N,\;T_j(\omega)\leq T,$ for almost every $\omega\in\Omega$.\\
The convergence $G_n\to G,\;n\to\infty,$ in $L^2([0,T]\times\mathcal O\times\R,\,dt\,d\xi\,\nu(d\sigma);\,H)$ implies the convergence $G_{n_j}(t,\xi,\sigma)\to G(t,\xi,\sigma),\;j\to\infty,$ in $H$ for $ds\,d\xi\,\nu(d\sigma)$-almost every $(t,\xi,\sigma)\in [0,T]\times\mathcal O\times\R$ (where $(G_{n_j})_{j\in\N}$ is a subsequence of $(G_n)_{n\in\N}$). This, the uniform convergence following from $\eqref{glmKgz}$ and the fact that the laws of the random vectors $(T_j,\Xi_j,\Sigma_j),\;j\in\N,$ are absolutely continuous with respect to $dt\,d\xi\,\nu(d\sigma)$ imply that, for almost every $\omega\in\Omega$, the jumps of the process $\big(\bar Y_{h,t}\big)_{t\in[0,T]}$ occur at the jump times $T_j(\omega),\;j\in\N,\;T_j(\omega)\leq T$ and are exactly $G\big(T_j(\omega),\Xi_j(\omega),\Sigma_j(\omega)\big)$, $j\in\N,\;T_j(\omega)\leq T.$\\
Consequently, we have
\begin{align*}
    \lefteqn{\E\sum_{t\leq T}\Big\{v_h(T-t,\bar Y_{h,t})-v_h(T-t,\bar Y_{h,t-}) -\left\langle D_xv_h(T-t,\bar Y_{h,t-}),\,\Delta\bar Y_{h,t}\right\rangle_H\Big\}}\notag\\
   &= \E\sum_{j\in\N,\,T_j\leq T}\Big\{v_h\big(T-T_j,\bar Y_{h,T_j-}+G(T_j,\Xi_j,\Sigma_j)\big)-v_h(T-T_j,\bar Y_{h,T_j-})\notag\\
    &\qquad\qquad\qquad\qquad\qquad-\left\langle D_xv_h(T-T_j,\bar Y_{h,T_j-}),\,G(T_j,\Xi_j,\Sigma_j)\right\rangle_H\Big\}.
\end{align*}
With similar arguments, one can show that
\begin{align*}
\lefteqn{\E\int_0^T\int_{\mathcal O\times \R}\Big\{v_h\big(T-t,\bar Y_{h,t-}+G(t,\xi,\sigma)\big)-v_h(T-t,\bar Y_{h,t-})}\\
&\qquad\qquad\qquad\qquad\qquad-\left\langle D_xv_h(T-t,\bar Y_{h,t-}),\,G(t,\xi,\sigma)\right\rangle_H\Big\}\,\pi(dt,d\xi,d\sigma)\\
&=\E\sum_{j\in\N,\,T_j\leq T}\Big\{v_h\big(T-T_j,\bar Y_{h,T_j-}+G(T_j,\Xi_j,\Sigma_j)\big)-v_h(T-T_j,\bar Y_{h,T_j-})\\
    &\qquad\qquad\qquad\qquad\qquad-\left\langle D_xv_h(T-T_j,\bar Y_{h,T_j-}),\,G(T_j,\Xi_j,\Sigma_j)\right\rangle_H\Big\}.
\end{align*}
This finishes the proof of $(\ref{jump-sum-trans})$.

\section{Extension to the case $\int_{\{|\sigma|\geq 1\}}\sigma^2\,\nu(d\sigma)=\infty$}\label{A2}

Here we show how to handle the case where assumption \eqref{AssNu0} fails to hold, assuming just that $\nu(\{0\})=0$ and
\begin{equation}\label{Levymeas}\tag{\ref{AssNu0}$'$}
\int_\R\min(\sigma^2,1)\,\nu(d\sigma)<\infty.
\end{equation}
$(Z_t)_{t\in[0,T]}$ as defined in Section 2 is a process with values in a Hilbert space $U$, which we assume to be of the form $U=H^{-\frac d2-\epsilon}(\mathcal O),\;\epsilon>0,$ in this section. Informally, $Z_t$ is also denoted by
\begin{align*}
Z(t,d\xi)&=\int_0^t\int_\R\sigma\,\hat\pi(ds,d\xi,d\sigma).
\end{align*}
As we cannot follow the construction in Section 2 if assumption \eqref{AssNu0} does not hold, we consider in this section a process of the form
\begin{align*}
Z(t,d\xi)&=\int_0^t\int_{\{|\sigma|\leq 1\}}\sigma\,\hat\pi(ds,d\xi,d\sigma)+\int_0^t\int_{\{|\sigma|> 1\}}\sigma\,\pi(ds,d\xi,d\sigma)\\
&=: M(t,d\xi)+ P(t,d\xi).
\end{align*}
More rigorously, we define the process $(Z_t)_{t\in[0,T]}$ by $Z_t:=M_t+P_t,\;t\in[0,T]$, where $(M_t)_{t\in[0,T]}$ is an impulsive cylindrical process with jump size intensity measure $\nu_M(d\sigma)=\one_{[-1,1]}(\sigma)\nu(d\sigma)$ and $(P_t)_{t\in[0,T]}$ is a $U$-valued compound Poisson process given by
\[P_t=\sum_{\substack{T_j\leq t\\ |\Sigma_j|>1}}\Sigma_j\cdot\delta_{\Xi_j}.\]
Here $\big((T_j,\Xi_j,\Sigma_j)\big)_{j\in\N}$ is the sequence of random elements in $[0,T]\times\mathcal O\times\R$ introduced in Section 2 and for $\xi\in\mathcal O$ we denote by $\delta_\xi\in U$ the Dirac delta function at $\xi$.

The equation
\begin{equation}\label{genEqn}\tag{\ref{sde}$'$}
dX_t+AX_t=Q^{1/2}\,dZ_t=Q^{1/2}\,(dM_t+dP_t),\quad X_0=x_0\in H,\quad t\in [0,T],
\end{equation}
has a weak solution with values in $H$ if the linear operator $Q^{1/2}:H\to H$ has a (not necessarily bounded) extension $Q^{1/2}:U\to H$, compare Example \ref{extendedOp} below. In this case a weak solution to \eqref{genEqn} can be constructed as follows (compare \cite{PesZab}, Chapter 9.7). Let us consider a sequence of stopping times
\[\tau_m:=\inf\{t\in[0,T]\,:\,P_t-P_{t-}\notin V_m\},\quad m\in\N, \]
where the $V_m$ are bounded measurable subsets of $U$ given by
\[V_m:=\left\{\sigma\delta_\xi\,:\,\xi\in\mathcal O,\,\sigma\in[-m,m]\right\},\quad m\in\N,\]
and $\inf\emptyset:=\infty$. Note that due to \eqref{Levymeas} we have $\wP\left(\bigcup_{m\in\N}\{\tau_m=\infty\}\right)=1$. For every $m\in\N$ we define
\[
P^{[m]}_t:=P_{\min(t,\tau_m)},\quad M^{[m]}_t:=P^{[m]}_t-tu_m,\quad t\in[0,T],
\]
with $u_m:=\E P^{[m]}_1\in U$, and consider the problem
\begin{equation}\label{stoppedEqn}\tag{\ref{sde}$''$}
dX^{[m]}_t+AX^{[m]}_t\,dt=Q^{1/2}\left(u_m\,dt+dM_t+dM_t^{[m]}\right),\quad X^{[m]}_0=x_0\in H.
\end{equation}
$\big(M^{[m]}_t\big)_{t\in[0,T]}$ is an impulsive cylindrical process in the sense of Section 2, with jump size intensity measure $\one_{[-m,m]\setminus[-1,1]}(\sigma)\,\nu(d\sigma)$.
There exists a unique weak solution $\big(X^{[m]}_t\big)_{t\in[0,T]}$ to \eqref{stoppedEqn} for every $m\in\N$, given by
\begin{equation*}
X^{[m]}_t=S(t)x_0+\int_0^tS(t-s)Q^{1/2}\,\left(u_m\,ds+dM_s+dM_s^{[m]}\right),\qquad t\in[0,T].
\end{equation*}
For $m\leq n\in\N$ one has $\{\tau_m=\infty\}\subset\{\tau_n=\infty\}$, and for every $t\in[0,T]$ the equality $X^{[m]}_t\one_{\{\tau_m=\infty\}}=X^{[n]}_t\one_{\{\tau_m=\infty\}}$ holds \wP-almost surely.
Consequently, the process $(X_t)_{t\in[0,T]}$ given by the $\wP$-almost surely unique limit
\[X_t:=\lim_{m\to\infty}X^{[m]}_t\one_{\{\tau_m=\infty\}}\]
is a weak solution to equation \eqref{genEqn}.

\begin{ex}\label{extendedOp}
Equation \eqref{genEqn} has an $H$-valued weak solution $(X_t)_{t\in[0,T]}$ if the operator $Q^{1/2}$ can be extended to $Q^{1/2}:U\to H$. Here we give an example for such an extension.
Consider an operator $Q$ as in the example mentioned in the introduction, \frenchspacing{i.e.}
\[Q^{1/2}x(\xi)=\int_\mathcal O q_0(\xi,\zeta)x(\zeta)\,d\zeta,\qquad x\in L^2(\mathcal O),\;\xi\in\mathcal O,\]
where $q_0$ is a positive semidefinite symmetric function on $\mathcal O\times\mathcal O$. Assume furthermore that $Q^{1/2}\in L(H,U^*)=L\big(L^2(\mathcal O),H^{d/2+\epsilon}_0(\mathcal O)\big)$. Then $Q^{1/2}$ can be extended to a linear and bounded operator $\in L(U,H)$. For $u\in U$ and $x\in H$ set
\[\langle Q^{1/2}u,x\rangle_H:={\phantom{\big(}}_{\scriptstyle U^*}\Big(u\,,\,\int_{\mathcal O}q_0(\,\cdot\,,\eta)x(\eta)\,d\eta\Big){\phantom{\big(}}_{\hspace{-5pt}U}={\phantom{\big)}}_{\scriptstyle U^*}\big(u,Q^{1/2}x\big)_U\,,\]
$\phantom{(}_{\scriptscriptstyle U^*}(\,\cdot\,,\,\cdot\,)_{\hspace{-1pt}\scriptscriptstyle U}:U^*\times U\to\R$ being the canonical dual pairing. By the Riesz representation theorem, this defines $Q^{1/2}:U\to H$. In particular, for $\xi\in\mathcal O$ we have the equality in $H$
\[Q^{1/2}\delta_\xi=q_0(\,\cdot\,,\xi).\]
\end{ex}

\begin{rem}
If the extension of $Q^{1/2}$ is continuous, i.e. if $Q^{1/2}\in L(U,H)$, then the term $Qu_m$ appearing in equation \eqref{stoppedEqn} can be represented as follows: Let $\mu_m$ denote the jump intensity measure of $\big( P^{[m]}_t\big)_{t\in[0,T]}$. We have
\begin{align*}
u_m =\int_Uz\,\mu_m(dz) &=\int_{\mathcal O\times([-m,m]\setminus[-1,1])}\sigma\delta_\xi\,d\xi\,\nu(d\sigma)\\
&=\int_{\mathcal O\times([-m,m]\setminus[-1,1])}\left\{\sum_{k=1}^\infty\sigma e_k(\xi)e_k\right\}\,d\xi\,\nu(d\sigma),
\end{align*}
where all integrals are $U$-valued $L^1$-Bochner integrals, $(e_k)_{k\in\N}$ is an orthonormal basis in $H$ and the infinite sum is a limit in the space $L^2\big(\mathcal O\times([-m,m]\setminus[-1,1]),d\xi\,\nu(d\sigma)\,;\,U\big)\subset L^1\big(\mathcal O\times([-m,m]\setminus[-1,1]),d\xi\,\nu(d\sigma)\,;\,U\big)$.\\
Consequently,
\[Q^{1/2}u_m=\sum_{k=1}^\infty\left\{\int_{\mathcal O\times([-m,m]\setminus[-1,1])}\sigma e_k(\xi)\,d\xi\,\nu(d\sigma)\right\}Q^{1/2}e_k,\]
the infinite sum being a limit in $H$.
\end{rem}

We extend the approximation scheme \eqref{scheme}--\eqref{theta} of Section 4 and define discretizations $(X^{n,[m]}_h)_{n\in\{0,\ldots,N\}}$ of $(X_t)_{t\in[0,T]}$ for integers $N,\;m\,\geq 1$ and real $h>0$ by
\begin{equation}\label{discrSoln'}\tag{\ref{discrSoln}$'$}
X_h^{n,[m]}:= S_{h,\Delta t}^n P_hx_0 +\sum_{k=0}^{n-1}S_{h,\Delta t}^{n-k-1}T_{h,\Delta t}P_h Q^{1/2}\Big(\Delta t\,u_m + (M_{t_{k+1}}-M_{t_k}) + \big(M^{[m]}_{t_{k+1}}-M^{[m]}_{t_k}\big)\Big),
\end{equation}
where $n\in\{0,\ldots,N\}$, compare \eqref{discrSoln}.
Given a bounded function $\varphi:H\to\R$ we find
\begin{align*}
\lefteqn{\left|\E\varphi\big(X_h^{N,[m]}\big)-\E\varphi\big(X_T\big)\right|}\\
&\qquad\leq \left|\E\varphi\big(X_h^{N,[m]}\big)-\E\varphi\big(X^{[m]}_T\big)\right| + \left|\E\varphi\big(X_T^{[m]}\big)-\E\varphi\big(X_T\big)\right|\\
&\qquad\leq \left|\E\varphi\big(X_h^{N,[m]}\big)-\E\varphi\big(X^{[m]}_T\big)\right| + 2\|\varphi\|_\infty \wP(\tau_m<\infty)\\
&\qquad\leq \left|\E\varphi\big(X_h^{N,[m]}\big)-\E\varphi\big(X^{[m]}_T\big)\right|+
2\|\varphi\|_\infty\min\big(T\,\lambda(\mathcal O)\,\nu(\R\setminus[-m,m])\,,\,1\big),
\end{align*}
where $\|\varphi\|_\infty:=\sup_{x\in H}|\varphi(x)|_H$ and $\lambda(\mathcal O)$ is the Lebesgue measure of $\mathcal O$. Of course, due to \eqref{Levymeas}, $\nu(\R\setminus[-m,m])$ tends to zero as $m\to\infty$ and the order of convergence depends on the particular jump size intensity measure $\nu$.

The term $\big|\E\varphi\big(X_h^{N,[m]}\big)-\E\varphi\big(X^{[m]}_T\big)\big|$ can be treated with the methods developed in Sections 4, 5 and 6. We briefly sketch the alterations. Let
\begin{align}
\tag{\ref{Xht}$'$}\label{Xht'}
X_{h,t}^{[m]} &= S_h(t)P_hx_0 + \int_0^tS_h(t-s)P_hQ^{1/2}\,\left(u_m\,ds+dM_s + dM^{[m]}_s\right),\qquad t\in [0,T],\\
\tag{\ref{Yht}$'$}
Y_{h,t}^{[m]} &:= S_h(T)P_hx_0+\int_0^t\Phi(s)\,\left(u_m\,ds+dM_s + dM^{[m]}_s\right),\qquad\; t\in [0,T],\\
\tag{\ref{Ybarht}$'$}
\bar{Y}_{h,t}^{[m]}&:= S^N_{h,\Delta t}P_hx_0 +\int_0^t\Gamma(s)\,\left(u_m\,ds+dM_s + dM^{[m]}_s\right),\qquad\; t\in [0,T],
\end{align}
with $\Phi(s)=S_h(T-s)P_hQ^{1/2}$ and
$\Gamma(s)=\sum_{k=0}^{N-1}S_{h,\Delta t}^{N-k-1}T_{h,\Delta t}\one_{(t_k,t_{k+1}]}(s)P_hQ^{1/2}$ as in Section 4. Let $v_h^{[m]}:[0,T]\times H\to\R$ and $v^{[m]}:[0,T]\times H\to\R$ be defined by
\begin{align}
v_h^{[m]}(t,x) &:= \E\varphi\left( x+ \int_{T-t}^T S_h(T-r)P_hQ^{1/2}\left(u_m\,ds+dM_s + dM^{[m]}_s\right)\right),\tag{\ref{defvh}$'$}\\
v^{[m]}(t,x) &:= \E\varphi\left( x+ \int_{T-t}^T S(T-r)Q^{1/2}\left(u_m\,ds+dM_s + dM^{[m]}_s\right)\right).\tag{\ref{defv}$'$}
\end{align}
Arguing along the lines of Section 5, we get
\begin{errorRep'}
Assume that \eqref{Levymeas} holds, $\|S(t)Q^{1/2}\|_{\text{(HS)}}\in L^2([0,T],dt)$ and that $Q^{1/2}$ has a linear and continuous extension $Q^{1/2}:U\to H$. Let $\varphi\in C^2_b(H)$ such that for each $x\in H$ the derivative $D^2\varphi(x)$ is an element of $L_{\text{(HS)}}(H)$ and that the mapping $x\mapsto D^2\varphi(x)\in L_{\text{(HS)}}(H)$ is uniformly continuous on any bounded subset of $H$. Let $T\geq 1,\;m\in\N$ and $(X_t^{[m]})_{t\in[0,T]}$ be the $H$-valued stochastic process satisfying equation \eqref{stoppedEqn}. For any $N\geq 1$ and $h>0$, let $(X_h^{n,[m]})_{n\in\{0,\ldots, N\}}$ be given by \eqref{discrSoln'} and let $(X_{h,t}^{[m]})_{t\in[0,T]}$ be as in \eqref{Xht'}.\\
Then, the following error expansions hold for the time discretization error
\begin{align*}
\lefteqn{\E\varphi\big(X^{N,[m]}_h\big)-\E\varphi\big(X^{[m]}_{h,T}\big)
\;=\;\left\{v_h^{[m]}(T,S_{h,\Delta t}^NP_hx_0)-v_h^{[m]}(T,S_h(T)P_hx_0)\right\}}\\
&\qquad+\E\int_0^T\int_{\mathcal O\times [-m,m]}\bigg\{v_h^{[m]}\Big(T-t,\bar Y_{h,t-}^{[m]}+G(t,\xi,\sigma)\Big)-v_h^{[m]}\Big(T-t,\bar Y_{h,t-}^{[m]}+F(t,\xi,\sigma)\Big)+\\
&\qquad\qquad\qquad\qquad\;+\one_{[-1,1]}(\sigma)\left\langle D_xv_h^{[m]}\Big(T-t,\bar Y_{h,t-}^{[m]}\big),\,F(t,\xi,\sigma)-G(t,\xi,\sigma)\right\rangle_H\bigg\}\,dt\,d\xi\,\nu(d\sigma),
\end{align*}
\begin{samepage}
and for the space discretization error
\begin{align*}
\lefteqn{\E\varphi\big(X^{[m]}_{h,T}\big)-\E\varphi\big(X_T^{[m]}\big)
\;=\;\left\{v^{[m]}(T,S_h(T)P_hx_0)-v^{[m]}(T,S(T)x_0)\right\}}\\
&\qquad+\E\int_0^T\int_{\mathcal O\times [-m,m]}\bigg\{v^{[m]}\Big(T-t, Y^{[m]}_{h,t-}+F(t,\xi,\sigma)\Big)-v^{[m]}\Big(T-t,Y^{[m]}_{h,t-}+E(t,\xi,\sigma)\Big)+\\
&\qquad\qquad\qquad\qquad+\one_{[-1,1]}(\sigma)\left\langle D_xv^{[m]}\Big(T-t, Y_{h,t-}^{[m]}\Big),\,E(t,\xi,\sigma)-F(t,\xi,\sigma)\right\rangle_H\bigg\}\,dt\,d\xi\,\nu(d\sigma).
\end{align*}
\end{samepage}
\end{errorRep'}

\begin{rem}
In some cases the assumption $\|S(t)Q^{1/2}\|_{\text{(HS)}}\in L^2([0,T],dt)$ in Theorem \ref{errorRep}$'$ is already fulfilled as a consequence of $Q^{1/2}$ having a linear and continuous extension $Q^{1/2}:U\to H$, because the latter is equivalent to $Q^{1/2}\in L(H,U^*)=L\big(L^2(\mathcal O),H^{d/2+\epsilon}_0(\mathcal O)\big)$. As an example, let $\mathcal O=(0,1)^d$ and $(A,D(A))=\big(-\Delta\,,\,H^1_0(\mathcal O)\cap H^2(\mathcal O)\big)$. Setting $\beta:=\frac d2+\epsilon$, we have $A^{\beta/2}Q^{1/2}\in L(H)$ due to $Q^{1/2}\in L(H,U^*)$ and this yields $\|S(t)Q^{1/2}\|_{\text{(HS)}}\in L^2([0,T],dt)$ because $\text{Tr}(A^{-\beta})<\infty$.
\end{rem}

If the jump size intensity measure $\nu$ satisfies
\begin{equation}\tag{\ref{AssNu}$'$}\label{AssNu'}
\int_\R\min(|\sigma|,1)\,\nu(d\sigma)<\infty,
\end{equation}
one gets the following counterpart of Theorem \ref{result}:
\begin{result'}
Assume \eqref{AssNu'} and \eqref{S1}, \eqref{S2}, \eqref{ass1}, \eqref{assTr}, \eqref{assQ} listed on page \pageref{result}. Assume furthermore that $Q^{1/2}$ has a linear and continuous extension $Q^{1/2}:U\to H$. Let $\varphi\in C^2_b(H)$ such that for each $x\in H$ the derivative $D^2\varphi(x)$ is an element of $L_{\text{(HS)}}(H)$ and that the mapping $x\mapsto D^2\varphi(x)\in L_{\text{(HS)}}(H)$ is uniformly continuous on every bounded subset of $H$. Let $T\geq 1$ and $(X_t)_{t\in[0,T]}$ be the $H$-valued stochastic process satisfying the equation
\begin{equation*}
dX_t+AX_t=Q^{1/2}\,dZ_t=Q^{1/2}\,(dM_t+dP_t),\quad X_0=x_0\in H,\quad t\in [0,T].
\end{equation*}
For any $N,\;m\geq 1$ and $h>0$, let $(X^{n,[m]}_h)_{n\in\{0,\ldots, N\}}$ be given by \eqref{discrSoln'}.\\
Then for all $\gamma<1-\alpha+\beta\leq 1$ there exists a constant $C=C(T,\varphi,A,Q,|x_0|,\gamma)>0$ which does not depend on $h,\;N$ and $m$ such that the following inequality holds
\begin{eqnarray*}
 \lefteqn{\big|\E\varphi\big(X^{N,[m]}_h\big)-\E\varphi(X_T)\big|}\\
 &\leq& C\int_{-m}^m|\sigma|\,\nu(d\sigma) (h^{2\gamma}+(\Delta t)^{\gamma})+
2\|\varphi\|_\infty\min\big(T\,\lambda(\mathcal O)\,\nu(\R\setminus[-m,m])\,,\,1\big),
\end{eqnarray*}
where $\Delta t =T/N\leq 1$.
\end{result'}
\end{appendix}
\bigskip
{\bf Acknowledgement:} Financial support by the Deutsche Forschungsgemeinschaft (grant SCHI-419/5-1) within the Priority Program 1324 is gratefully acknowledged.

{\frenchspacing

}
\end{document}